\documentclass{article}
\usepackage[cp1251]{inputenc}
\usepackage[english]{babel}
\usepackage[12pt]{extsizes}
\usepackage{mathrsfs}
\usepackage{amsfonts,amssymb}
\usepackage[tbtags]{amsmath}
\usepackage[T2A]{fontenc}
\usepackage[unicode, pdftex]{hyperref}
\usepackage{mathtools}
\usepackage{enumitem}
\usepackage{graphicx}
\usepackage{caption} 
\usepackage{subcaption}
\usepackage{float,wrapfig}
\usepackage{blindtext,titlefoot}
\DeclarePairedDelimiter{\floor}{\lfloor}{\rfloor}

\DeclarePairedDelimiter{\ceil}{\lceil}{\rceil}

\newtheorem{theorem}{Theorem}

\newtheorem{lemma}{Lemma}
\newtheorem{claim}{Claim}

\begin{document}

\title{Cycle saturation in random graphs}
\author{Yu. Demidovich\footnote{Moscow Institute of Physics and Technology (National Research University), 9 Institutskiy per., Dolgoprodny, Moscow Region, Russian Federation, e-mail: demidovich.yua@phystech.edu}\and A. Skorkin\footnote{Adyghe State University, ul. Pervomayskaya, 208, Maykop, Republic of Adygea, Russian Federation}\and M. Zhukovskii\footnote{Moscow Institute of Physics and Technology (National Research University), 9 Institutskiy per., Dolgoprodny, Moscow Region, Russian Federation; Adyghe State University, Caucasus mathematical center, ul. Pervomayskaya, 208, Maykop, Republic of Adygea, Russian Federation; The Russian Presidential Academy of National Economy and Public Administration, Prospect Vernadskogo, 84, bldg 2, Moscow, Russian Federation}}

\date{}
\maketitle\unmarkedfntext{\textit{In the following Eurocomb'21 Proceedings paper weaker results with no proofs were published:} \\Demidovich Y., Zhukovskii M. (2021) Cycle Saturation in Random Graphs. In: Ne\v{s}et\v{r}il J., Perarnau G., Ru\'{e} J., Serra O. (eds) Extended Abstracts EuroComb 2021. Trends in Mathematics, vol 14. Birkh\"{a}user, Cham., pp 811-816.}
\begin{abstract}
	For a fixed graph $F,$ the minimum number of edges in an edge-maximal $F$-free subgraph of $G$ is called the $F$-saturation number. The asymptotics of the $F$-saturation number of the binomial random graph $G(n,p)$ for constant $p\in(0,1)$ is known for complete graphs $F=K_m$ and stars $F=K_{1,m}.$ This paper is devoted to the case when the pattern graph $F$ is a simple cycle $C_m.$ We prove that, for $m\geqslant 5,$ whp $\mathrm{sat}\left(G\left(n,p\right),C_m\right) = n+\Theta\left(\frac{n}{\ln n}\right).$ Also we find $c=c(p)$ such that whp $\frac{3}{2}n(1+o(1))\leqslant\mathrm{sat}\left(G\left(n,p\right),C_4\right)\leqslant cn(1+o(1)).$ In particular, whp $\mathrm{sat}\left(G\left(n,\frac{1}{2}\right),C_4\right)\leqslant\frac{27}{14}n(1+o(1)).$
\end{abstract}
\section{Introduction}
Let $F$ be a graph. Zykov \cite{Zyk} and later independently Erd\H{o}s, Hajnal and Moon \cite{ErHajMoon} raised a question of finding the minimum number of edges in an edge-maximal $F$-free graph on $n$ vertices. Formally, a graph $H$ is said to be $F$-\textit{saturated} if it is a maximal $F$-free graph, i.e. $H$ does not contain any copy of $F$ as a subgraph, but adding any missing edge to $H$ creates one. The \textit{saturation number} $\mathrm{sat}(n,F)$ is defined to be the minimum number of edges in an $F$-saturated graph on $n$ vertices. 

If $F$ is an $m$-clique then $\mathrm{sat}(n,F)$ is known. It was proven in \cite{ErHajMoon} that when $n\geqslant m \geqslant 2,$ then
\begin{equation}\label{classic_complete}
	\mathrm{sat}(n,K_m) = (m-2)n - \binom{m-1}{2}.
\end{equation}

For $K_{1,m},$ the star graph on $m+1$ vertices, the saturation number is also known. It was proven by K\'{a}szonyi and Tuza \cite{KasTu} that
\begin{equation*}
	\mathrm{sat}(n,K_{1,m})=
	\begin{cases}
		\binom{m}{2}+\binom{n-m}{2}, & m+1\leqslant n\leqslant\frac{3m}{2};\\
		\ceil[\big]{\frac{(m-1)n}{2}-\frac{m^2}{8}},& n\geqslant\frac{3m}{2}.
	\end{cases}
\end{equation*}
Finding $\mathrm{sat}(n,C_m)$ is harder (as usual, $C_m$ is a simple cycle on $m$ vertices). The problem is completely solved only for $m=4, 5.$ It was determined in \cite{Oll} by Ollman that, for $n\geqslant 5,$
\begin{equation}\label{classic_cycle_4}
	\mathrm{sat}(n,C_4) = \floor[\bigg]{\frac{3n-5}{2}}.
\end{equation}
If $n\geqslant 21,$ then (\cite{Chen}, \cite{FFL})
\begin{equation}\label{classic_cycle_5}
	\mathrm{sat}(n,C_5)=\ceil[\bigg]{\frac{10}{7}(n-1)}.
\end{equation}

Luo, Shigeno and Zhang in \cite{LSZ} established that
\begin{equation}\label{classic_cycle_6}
	\mathrm{sat}(n,C_6)\leqslant\floor[\bigg]{\frac{3n-3}{2}},\ \ \ {\rm for}\ n\geqslant 9,\qquad
	\mathrm{sat}(n,C_6)\geqslant \ceil[\bigg]{\frac{7n}{6}}-2,\ \ \ {\rm for}\ n\geqslant 6.
\end{equation}

Finally, F\"{u}redi and Kim \cite{FK} showed, that for all $m\geqslant 7$ and $n\geqslant 2m-5,$
\begin{equation}\label{classic_cycle_m}
	\left(1+\frac{1}{m+2}\right)n-1<\mathrm{sat}(n,C_m)<\left(1+\frac{1}{m-4}\right)n+\binom{m-4}{2}.
\end{equation}

More results concerning the saturation problem can be found, e.g., in \cite{FFS} and in references therein.\\

Kor\'{a}ndi and Sudakov \cite{KorSud} initiated the study of the saturation problem for random graphs.

Recall that the \textit{random graph} $G(n,p)$ is a random element of the set of all graphs $G$ on $[n]:=\{1,\ldots,n\}$ with probability distribution ${\sf P}(G(n,p)=G)=p^{|E(G)|}(1-p)^{\binom{n}{2}-|E(G)|}$ (or, in other words, every pair of vertices is adjacent with probability $0\leqslant p \leqslant 1$ independently). We say that a graph property $Q$ holds with high probability (whp), if ${\sf P}\left(G(n,p)\in Q \right) \to 1$ as $n\to\infty.$

For fixed graphs $F$ and $G,$ we say that a spanning subgraph $H\subseteq G$ is $F$-\textit{saturated} in $G$ if $H$ is an inclusion-maximal $F$-free spanning subgraph of $G.$ The minimum number of edges in an $F$-saturated graph in $G$ is denoted by $\mathrm{sat}(G,F)$ (in particular, $\mathrm{sat}(n,F)=\mathrm{sat}(K_n,F)).$

Kor\'{a}ndi and Sudakov \cite{KorSud} asked a question of determining the saturation number of $G(n,p)$ when $F=K_m.$ They proved that, for every fixed $p\in(0,1)$ and fixed integer $m\geqslant 3,$ whp
\begin{equation}\label{complete}
	\mathrm{sat}(G(n,p),K_m)=(1+o(1))n\log_{\frac{1}{1-p}}n.
\end{equation}

The saturation number of $G(n,p)$ when $F$ is a star graph was studied in a couple of papers. Note that, by the definition, $\mathrm{sat}(G,K_{1,2})$ coincides with the minimum cardinality of a maximal matching in $G.$ Zito \cite{Zito} showed that whp
$\frac{n}{2}-\log_{\frac{1}{1-p}}(np)<\mathrm{sat}(G(n,p),K_{1,2})<\frac{n}{2}-\log_{\frac{1}{1-p}}\sqrt{n}.$

Notice that this result can be easily improved. If we fix a bipartition of $[n],$ it follows from a well-known bound for the probability of non-containing a perfect matching (see Theorem \ref{bip} in Section \ref{prelim}) that 
\begin{equation}\label{er}
	{\sf P}\left(G(n,p)\text{ has no perfect matching}\right) = O\left(ne^{-np/2}\right).
\end{equation}
Therefore, by the union bound, whp any induced subgraph of $G(n,p)$ on at least $n-2\log_{\frac{1}{1-p}}n$ vertices contains a perfect matching. From Lemma \ref{lemma_lower_3} (see Section \ref{prelim}), we can conlude that whp the induced subgraph of $G(n,p)$ obtained by removing a maximum independent set contains a perfect matching. Therefore, whp $\mathrm{sat}\left(G(n,p), K_{1,2}\right) = \ceil{(n-\alpha(G(n,p)))/2},$ i.e. is also concentrated in a set consisting of two consecutive values.

Mohammadian and Tayfeh--Rezaie \cite{MohTay} proved, that for every fixed $p\in(0,1)$ and fixed integer $m\geqslant 3,$ whp
\begin{equation*}
	\mathrm{sat}(G(n,p),K_{1,m})=\frac{(m-1)n}{2}-(1+o(1))(m-1)\log_{\frac{1}{1-p}}n.
\end{equation*}

When $F$ is a complete graph, the comparison of \eqref{classic_complete} and \eqref{complete} shows that the saturation number becomes roughly logarithm times bigger after the random deletion of edges. When $F$ is a star on $m+1$ vertices, there is an asymptotical stability of the saturation number.\\

It is natural to ask a question about an asymptotical behavior of the $C_m$-saturation number of $G(n,p).$ The first result of the present paper establishes an asymptotical behavior of the $C_m$-saturation number of $G(n,p)$ when $m\geqslant 5.$

\begin{theorem}\label{th_main_1}
	Let $p\in(0,1)$ be fixed. For every $m\geqslant 5,$ whp
	\begin{equation}\label{cycle_m5}
		n+\frac{n}{4(m-1)\log_{\frac{1}{1-p}} n}(1-o(1))\leqslant\mathrm{sat}\left(G\left(n,p\right),C_m\right)\leqslant n + \frac{n}{2\log_{\frac{1}{1-p}} n}(1+o(1)).
	\end{equation}
\end{theorem}

\noindent\textbf{Remark.} The proof of the upper bound is constructive. It is based on the fact that whp almost all vertices of $G(n,p)$ can be covered by induced $K_{1,a}$ with $a=2(1+o(1))\log_{\frac{1}{1-p}}n.$ It is possible to improve the second summand in the upper bound up to the $\frac{m-3}{m-2}$ factor by replacing $K_{1,a}$ with, so-called, \textit{sparkler graphs}. We do not give a proof of the improvement since it does not give any additional insight, but the computations are much more dense. Let us just discuss the main ingredient of the construction. Let $a=2(1+o(1))\log_{\frac{1}{1-p}}n$  where $o(1)$ is chosen in an appropriate way. $S_{b,a}$ is a sparkler graph which is obtained from $P_b$ and $K_{1,a}$ by identifying the central vertex of $K_{1,a}$ and an end vertex of $P_b,$ the other end vertex of $P_b$ is the root of the sparkler graph. The improvement of the upper bound is based on the fact that, whp, in $G(n,p),$ there exists a set $\left\lbrace S_1,\ldots,S_{(m-2)t}\right\rbrace,$ where $t = \frac{n(1+o(1))}{2(m-2)\log_{1/(1-p)}n},$ of induced and vertex-disjoint copies of $S_{\floor[\big]{\frac{m-1}{2}}-1,a}$ with roots adjacent to $1$ such that, for every $i\in[t],$ the ends of $S_{(i-1)(m-2)+1},\ldots,S_{i(m-2)}$ (denoted by $r_{(i-1)(m-2)+1},\ldots,r_{i(m-2)})$ are adjacent sequentially: $r_j\sim r_{j+1},$ $j\in\{(i-1)(m-2)+1,\ldots,i(m-2)-1\},$ and there are no other edges between them. Moreover, there are no edges between distinct $P_{\floor[\big]{\frac{m-1}{2}}-1}$ that belong to $S_{(i-1)(m-2)+1},\ldots,S_{i(m-2)}$ other than those in $N(1).$\\

A construction that gives the upper bound in Theorem \ref{th_main_1} is obtained from a graph of size $\Theta(n/\ln n)$ with minimum degree bigger than $1$ and a small diameter by attaching to its vertices disjoint stars (see Figure \ref{fig:graph}). This construction is, in some sense, optimal: in our proof of the lower bound, we show that, after the recursive deletion of vertices with degree $1$ from a $C_m$-saturated subgraph of $G(n,p)$, the final graph has $\Omega(n/\ln n)$ vertices, minimum degree bigger than $1$ and a bounded diameter. Clearly, the last argument immediately implies the lower bound $n+\Theta(n/\ln n)$.

So, in contrast to $K_n,$ for $G(n,p),$ we know the exact asymptotics of the $C_m$-saturation number for all $m\geqslant 5:$ whp
$\mathrm{sat}\left(G(n,p),C_m\right)=n+\Theta\left(\frac{n}{\ln n}\right).$

The comparison of \eqref{classic_cycle_5},\eqref{classic_cycle_6},\eqref{classic_cycle_m} and our result \eqref{cycle_m5} implies that the order of growth of the saturation number is stable for $m>4,$ but there is no asymptotical stability since the constants in front of $n$ are different. Remarkably, our results also demonstrate that the saturation number drops after the random deletion of edges when $F=C_m,$ $m\geqslant 5,$ which is not the case for $F=K_m$ or $F=K_{1,m}.$\vspace{0.2cm}

The second result of the paper provides an upper bound for the $C_4$-saturation number of $G(n,p).$
\begin{theorem}\label{th_main_2}
	Let $p\in(0,1)$ be fixed. Whp 
	\begin{equation}\label{cycle_4_up}
		\mathrm{sat}\left(G\left(n,p\right),C_4\right)\leqslant \frac{3(1+(1-p)^3)}{2(1-(1-p)^3)}n(1+o(1)),
	\end{equation}
	when $p>1-1/\sqrt[3]{7},$ and
	\begin{equation}\label{cycle_4'_up}
		\mathrm{sat}\left(G\left(n,p\right),C_4\right)\leqslant \left(\frac{s+1}{2}+s(1-p)^s+\frac{s(1-p)^{2s}}{1-(1-p)^s}\right)n(1+o(1)),
	\end{equation}
	when $p\leqslant1-1/\sqrt[3]{7},$ where $s$ is the minimum positive integer such that $(2s^2+1)(1-p)^s < 1.$
\end{theorem}
In particular, when $p=1/2,$ then, whp $\mathrm{sat}\left(G\left(n,p\right),C_4\right)\leqslant \frac{27}{14}n(1+o(1)).$\vspace{0.2cm}

The proofs of bounds \eqref{cycle_4_up} and \eqref{cycle_4'_up} are constructive. Both constructions are recursive. For every $i$, we divide the current set of vertices $V_i$ into two parts $V_{i+1}$ and $V_i\setminus V_{i+1}$ and describe those edges that have at least one end-point outside $V_{i+1}$ (these edges are presented in Figures \ref{fig:graphc4} and \ref{fig:graph_prop_M}), and then move to the set $V_{i+1}$ at the next step. The first set $V_1$ is the set of all vertices $[n]$, and the final set has cardinality $o(\sqrt{n})$ (and thus the number of edges in this set does not affect the bound). While the existence of such a subgraph in $G(n,p)$ when $p>1-1/\sqrt[3]{7}$ is more or less straightforward, for $p\leq 1-1/\sqrt[3]{7}$, the proof of the existence relies on a nice fact about the existence of a perfect matching in the bipartite binomial random graph with some restrictions on possible edges inside the matching (see Lemma \ref{Hall}).\\

The third result of the paper provides a lower bound for the $C_4$-saturation number of $G(n,p).$
\begin{theorem}\label{th_main_3}
	Let $p\in(0,1)$ be fixed. Then whp
	\begin{equation}\label{cycle_4_low}
		\mathrm{sat}\left(G\left(n,p\right),C_4\right)\geqslant\frac{3}{2} n(1+o(1)).
	\end{equation}
\end{theorem}

Despite the fact that the asymptotics of the bound in Theorem 3 equals the $C_4$-saturation number for $K_n,$ the proof of Theorem 3 is much harder. In particular, if $F$ is a $C_4$-saturated subgraph of $G(n,p)$ (but not $K_n$) and vertices $u,\,v$ are not adjacent in $F$ and have degrees at most $2,$ then we can not guarantee that there they are at distance at most $3$ in $F$ and thus the arguments that imply the value of $\mathrm{sat}(K_n,C_4)$ completely fail. Our proof follows from several bounds on the number of vertices with degrees at most $1$ and the number of induced paths in $C_4$-saturated graphs in certain subgraphs of the binomial random graphs (see Lemma \ref{lemma_c4_2} and Lemma \ref{lemma_c4_3}) as well as some facts about the structure of the neighborhood $W$ of the set of vertices $U$ with degree $2$ and both neighbors with degrees bigger than $2.$ The crucial fact about the sets $U$ and $W$ is that whp the set of vertices from $U$ such that their neighbors have bounded degrees has size $o(n).$

Analysis of \eqref{classic_cycle_4} and \eqref{cycle_4_up},\eqref{cycle_4'_up},\eqref{cycle_4_low} shows that we can not make a conclusion about the stability of the $C_4$-saturation number, but we can infer that, after the random deletion of edges, the saturation number does not diminish asymptotically in comparison to the case of $F=C_m,$ $m\geqslant 5.$\\

Here we outline the structure of our work. In Section \ref{prelim} several auxiliary results used throughout the paper are provided. In Section \ref{proof_th2} we prove Theorem \ref{th_main_1}. The proof of Theorem \ref{th_main_2} is given in Section \ref{proof_th3}. Our last result, Theorem \ref{th_main_3}, is proven in Section \ref{proof_th4}.

\section{Preliminaries}\label{prelim}
Recall that the \textit{random graph} $G(n,n,p)$ is a random graph obtained from the complete bipartite graph $K_{n,n}$ by independent removal of each edge with probability $1-p.$ Below we frequently use the following result on the existence of a perfect matching in $G(n,n,p)$ (see \cite[Remark 4.3]{JLR}).
\begin{theorem}\label{bip}
	${\sf P}\left(G(n,n,p)\text{ has no perfect matching}\right)=O(ne^{-np}).$
\end{theorem}

Also we make use of the following result on the concentration of the independence number of $G(n,p).$
\begin{lemma}\label{lemma_lower_3}
	Let $p\in(0,1)$ be constant. The independence number of $G(n,p)$ whp belongs to $\{f(n),f(n)+1\},$ where 
	$$
	f(n)=\floor[\bigg]{ 2\log_{\frac{1}{1-p}}n-2\log_{\frac{1}{1-p}}\log_{\frac{1}{1-p}}n+2\log_{\frac{1}{1-p}}(e/2)+0.9}.
	$$
\end{lemma}  
For the proof we refer the reader to \cite[Theorem 11.1]{BolBook}.

We use the following version of Chernoff's bound (see \cite[Theorem 2.1]{JLR}).
\begin{theorem}\label{Chernoff} Suppose that $X$ is a binomial random variable with parameters $n,$ $p.$ Put $\lambda=np.$ Then, for every $t\geqslant 0,$
$$
{\sf P}\left(X\geqslant \lambda+t\right)\leqslant\exp\left(-\frac{t^2}{2(\lambda+t/3)}\right),
$$
$$
{\sf P}\left(X\leqslant \lambda-t\right)\leqslant\exp\left(-\frac{t^2}{2\lambda}\right).
$$	
\end{theorem}

Let $\Gamma$ be a graph on the vertex set $[n].$ Denote by $N_{\Gamma}(v_1,\ldots,v_{\ell})$ the set of all common neighbors of $v_1,\ldots,v_{\ell},$ $\ell\in\mathbb{N},$ in $\Gamma.$ Define $\overline{\Gamma}$ as the complement graph of $\Gamma,$ i.e., the graph on $[n]$ where every edge appears iff it is not present in $\Gamma.$ Notice that $N_{\overline{\Gamma}}(v_1,\ldots,v_{\ell})$ is the set of all common non-neighbors of $v_1,\ldots,v_{\ell}.$ Let $\Gamma[V]$ stand for the graph induced by $\Gamma$ on the vertex set $V\subset [n],$ i.e. the graph on $V,$ where every edge appears iff it is present in $\Gamma$ on $V.$ 
\section{Proof of Theorem \ref{th_main_1}}\label{proof_th2}
\subsection{Lower bound}

We start with the formulation of several helpful properties of $G(n,p).$
\begin{lemma}\label{lemma_lower_1}
	There exists a constant $A>0$ such that, for $p>A\sqrt{\frac{\ln n}{n}},$
	whp every pair of vertices has a common neighbor in $G(n,p).$
\end{lemma}
Its proof can be found in \cite[Theorem 2]{Spencer}.$\\$

\begin{lemma}\label{lemma_lower_2}
	Let $p\in(0,1)$ be a constant. There exists a constant $C>0$ such that whp the maximum size of an induced tree in $G(n,p)$ belongs to $\{f_C(n)-1,f_C(n)\},$ where
	$$
	f_C(n)=\floor{2\log_{\frac{1}{1-p}}n+C}.
	$$
\end{lemma}
It is proven in \cite[Theorem 2]{KSZ}.

\begin{lemma}\label{lemma_lower_4}
	There exists $\delta>0$ such that whp, after deletion of any $\floor{\delta n}$ vertices from $G(n,p),$ its diameter does not change.
\end{lemma}

Its proof easily follows from the following observation. By Theorem \ref{Chernoff} and by the union bound over all possible pairs of vertices, whp any two vertices have $np^2(1+o(1))$ common neighbors. Therefore, we can choose any $0<\delta<p^2.$\\

Let $\varepsilon > 0,$ $G$ be a graph on $[n],$ $n$ is large enough. Assume that $G$ has the property defined in Lemma \ref{lemma_lower_3}, the property defined in Lemma \ref{lemma_lower_1}, the property defined in Lemma \ref{lemma_lower_2} (for some $C>0$), the property defined in Lemma \ref{lemma_lower_4} (for some $\delta>0$). It remains to prove that
\begin{equation}\label{eq_sat}
	\mathrm{sat}\left(G,C_m\right)\geqslant n+\frac{n}{4(m-1)\log_{1/(1-p)}n}(1-\varepsilon).
\end{equation}
Fix a $C_m$-saturated graph $H$ in $G.$ First, we show that its diameter $d$ is at most $2m-2.$ Indeed, suppose the opposite. We know that any two vertices in $G$ have a common neighbor. Consider two vertices at distance $d$ in $H.$ Since $H$ is $C_m$-saturated in $G,$ their common neighbor in $G$ must be at  distance at most $m-1$ from each of them in $H.$ Hence, there exists a shorter path between them, which leads to a contradiction.

Starting from the graph $H,$ iteratively delete all vertices of degree $1$ until the subgraph of $H$ with minimum degree at least $2$ remains. Denote the final graph by $\tilde{H}.$ Observe that, for each deleted vertex $v,$ there exists a unique inclusion-maximal subtree $T(v)$ such that $v\in T(v)$ and $T(v)$ shares exactly $1$ vertex with $\tilde{H}$ (below, we call this vertex \textit{the root of }$T(v)$). Clearly, the diameter of $\tilde{H}$ is still at most $2m-2.$ 

Let us prove several properties of deleted trees.

\begin{claim}\label{claim1}
	Every vertex of a tree $T(v)$ can not have more than $\floor{2\log_{1/(1-p)}n}$ children.
\end{claim}
\textit{Proof.} Assume the contrary: suppose there exists a tree with a vertex that has more than $\floor{2\log_{1/(1-p)}n}$ children. By the respective property, such a set has at least one edge in $G$ which leads to a contradiction since there is no $P_m$ in $H$ connecting both children.
\begin{flushright}
	$\blacksquare$
\end{flushright}

The \textit{height} of a tree is the length (the number of edges) of the longest path from its root $R.$ Let $h=h(T)$ denote the height of $T=T(v).$ Let $L_i(T),$ $i\in[h],$ be the set of the vertices at distance $i$ from the root of $T.$

\begin{claim}\label{claim2}
	For every removed vertex $v,$ we have $h(T(v))\leqslant m-2.$
\end{claim}
\textit{Proof.} Suppose there exists a removed vertex $v$ with $h(T(v))>m-2.$ Let $X=X(T)=L_{m-1}(T)\sqcup\ldots\sqcup L_{h}(T).$ Consider the set $W$ of vertices from $[n]\setminus X$ adjacent to at least $1$ vertex from $X$ in $G.$ By the respective property, there exists $\delta>0$ such that after deletion of any $\floor{\delta n}$ vertices the diameter of $G$ remains equal to $2,$ which means that it remains connected. Hence, $|W|>\floor{\delta n}.$ But $W\subset \{R\}\cup L_1(T)\sqcup\ldots\sqcup L_{m-2}(T)$ since any vertex outside $T(v)$ is at distance at least $m$ from every vertex of $X.$ Then the tree $T[\{R\}\cup L_1(T)\cup\ldots\cup L_{m-2}(T)]$ of height $m-2$ has a vertex with more than $\floor{2\log_{1/(1-p)}n}$ children and we come into a contradiction with Claim \ref{claim1}.
\begin{flushright}
	$\blacksquare$
\end{flushright}
\begin{claim}\label{claim3}
	There exist at most $\floor{2\log_{1/(1-p)}n}$ distinct $T(v)$ such that $h(T(v))>\frac{m-2}{2}.$
\end{claim}
\textit{Proof.} Assume the contrary: suppose that there exist more than $\floor{2\log_{1/(1-p)}n}$ such trees. In every such tree, include exactly one vertex at distance $\floor{\frac{m-2}{2}}+1$ from the root into a set $A.$ By the condition on $\alpha\left(G\right),$ $A$ has two vertices connected by an edge which leads to a contradiction since there is no $P_m$ between them in $H.$
\begin{flushright}
	$\blacksquare$
\end{flushright}
\begin{claim}\label{claim4}
	For any deleted $v$ such that $h(T(v))\leqslant \frac{m-2}{2},$ we have $|V(T(v))|\leqslant f_C(n).$
\end{claim}
\textit{Proof.} If a tree $T(v)$ has height at most $\frac{m-2}{2},$ then it is induced by $G$ since any pair of its vertices is at distance at most $m-2$ in $H.$ Due to the respective property, there are no induced trees of size $f_C(n)+1$ in $G.$
\begin{flushright}
	$\blacksquare$
\end{flushright}

From Claims \ref{claim2}, \ref{claim3} and \ref{claim4}, we get that, for all but at most $\left(f_C(n)\right)^{m-1}$ vertices $v,$ $|V(T(v))|\leqslant f_C(n).$

Therefore, 
$$
|V(\tilde{H})|=:x\geqslant \frac{n-\left(f_C(n)\right)^{m-1}}{f_C(n)}>\frac{n(1-\varepsilon)}{2\log_{1/(1-p)} n}+1.
$$ A graph on $x$ vertices with minimum degree at least $2$ and diameter at most $2m-2$ has at least $a = x-1+\frac{x-1}{2(m-1)}$ edges. Indeed, a spanning tree of such a graph with diameter at most $2m-2$ has $x-1$ edges. In the worst case, all paths from a central vertex of $T$ to the leaves do not branch. Hence, the number of leaves is at least $b =\frac{x-1}{m-1}.$ In turn, at least $b/2$ edges are required to eliminate $b$ leaves.

Therefore, 
$$
|E(H)|\geqslant n+\frac{n}{4(m-1)\log_{1/(1-p)} n}(1-\varepsilon).
$$
Inequality \eqref{eq_sat} follows.

\subsection{Upper bound}
Consider a graph $\Gamma$ and a graph $H.$ An \textit{induced} $H$\textit{-factor} in $\Gamma$ is a set of induced and vertex-disjoint isomorphic copies of $H$ that cover all the vertices of $\Gamma.$ Let $v\in V\left(\Gamma\right).$ Let us call a set $\{S_1,\ldots,S_{\ell}\}$ of induced and vertex-disjoint copies of $K_{1,a}$ (different stars may have different sizes, and $a$ may be equal to $0$) in $\Gamma$ a $v$-\textit{factor}, if 
\begin{itemize}
	\item central vertices of $S_1,\ldots,S_{\ell}$ are adjacent to $v;$
	\item $V(S_1\sqcup \ldots\sqcup S_{\ell})=[n]\setminus \{v\}.$
\end{itemize}

\begin{lemma}\label{upper}
	Whp, in $G(n,p),$ there exists a $1$-factor $\left\lbrace S_1,\ldots,S_{\ell}\right\rbrace$ such that $\ell = \frac{n}{2\log_{1/(1-p)}n}(1+o(1))$ and the set of central vertices of $S_1,\ldots,S_\ell$ induces a subgraph that contains an induced $C_{m-2}$-factor.
\end{lemma}

\begin{figure}
	\centering
	\includegraphics[width=0.4\linewidth]{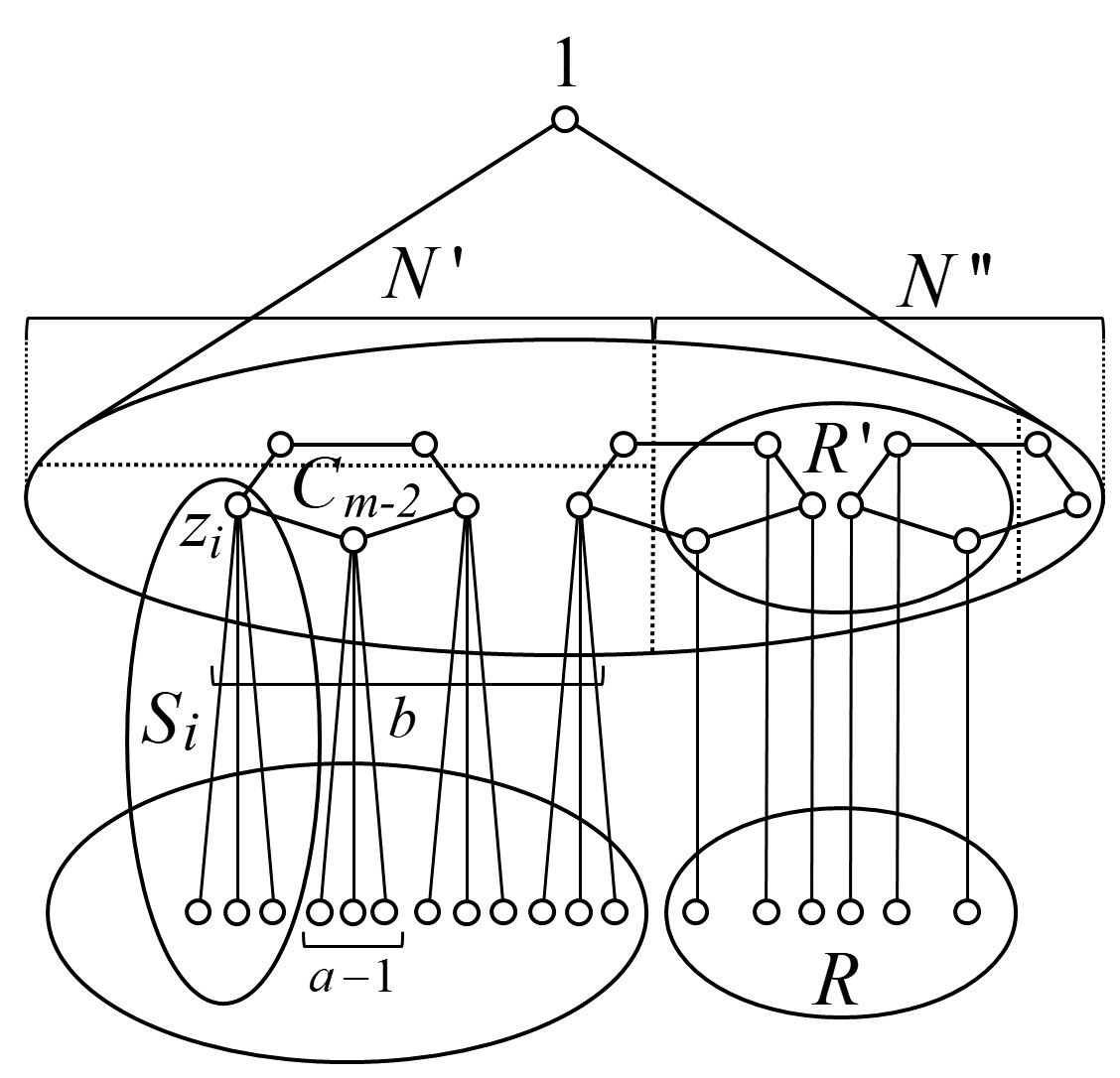}
	\caption{}
	\label{fig:graph}
\end{figure}

Let us first finish the proof of Theorem \ref{th_main_1} and then switch to the proof Lemma \ref{upper}. Clearly, the subgraph of $G(n,p)$ with edges
\begin{itemize}
	\item between $1$ and the central vertices of $S_1,\ldots,S_{\ell},$
	\item of an induced $C_{m-2}$-factor in the subgraph induced by the set of central vertices of $S_1,\ldots,S_{\ell},$
	\item of stars $S_1,\ldots,S_{\ell}$
\end{itemize}
(see Figure \ref{fig:graph}) is $C_m$-saturated in $G(n,p)$ and has exactly $n+\frac{n}{2\log_{1/(1-p)}n}\left(1+o(1)\right)$ edges.
This finishes the proof of the upper bound in Theorem \ref{th_main_1}.\\

\noindent\textbf{Proof of Lemma \ref{upper}.} Let $a=\floor{2\log_{1/(1-p)}n-8\log_{1/(1-p)}\ln n}.$ Let $b$ be the minimum integer such that $ab+3\floor[\big]{\frac{n}{\ln^2n}}+2m>n.$ Let $\ell\leqslant b+2\floor[\big]{\frac{n}{\ln^2 n}}+2m-1$ be the maximum integer divisible by $m-2.$
Expose all edges of $G(n,p)$ adjacent to $1.$ Choose an arbitrary set $N\subset N(1)$ of size $\ell.$
\begin{claim}\label{claim6}
	Whp $G(n,p)[N]$ has an induced $C_{m-2}$-factor.
\end{claim}

\noindent\textbf{Remark.} We suggest that Claim \ref{claim6} is known and its proof has already appeared somewhere. Unfortunately, we have not found it, and, by this reason, give its proof below. Notice that the notion of induced factors differs from the notion of (not necessarily induced) factors. The problem of existence of the latter is well-studied, (see, e.g. \cite{AY}, \cite{JKV}, \cite{Ruc}).

\begin{claim}\label{claim5}
	Whp, for every disjoint $A\subset N(1),$ $B\subset[n]\setminus\{1\}$ of size $|A|=|B|=\floor[\big]{\frac{n}{\ln^2n}}$ there exists an induced $K_{1,a-1}$ with a central vertex in $A$ and all the other vertices in $B.$
\end{claim}
	
Let us first finish the proof of Lemma \ref{upper} and, after that, proceed with the proofs of Claims \ref{claim6} and \ref{claim5}.

Find induced cycles $U_1,\ldots,U_t$ as stated in Claim \ref{claim6}. Let $N = N'\sqcup N'',$ where $|N''| = \ell - (b+\floor[\big]{\frac{n}{\ln^2n}}-1).$ Choose $V_1\subset N'$ of size $\floor[\big]{\frac{n}{\ln^2n}}$ and $W_1\subset [n]\setminus\left(\{1\}\cup N\right)$ of size $\floor[\big]{\frac{n}{\ln^2n}}.$ Find $S_1= K_{1,a-1}$ (with a central vertex $z_1$) in $V_1\sqcup W_1$ as stated in Claim \ref{claim5}. Choose $V_2\subset N'\setminus\{z_1\}$ of size $\floor[\big]{\frac{n}{\ln^2n}}$ and $W_2\subset [n]\setminus\left(\{1\}\cup N\cup V(S_1)\right)$ of size $\floor[\big]{\frac{n}{\ln^2 n}}.$ Find $S_2=K_{1,a-1}$ (with a central vertex $z_2$) in $V_2\sqcup W_2$ as stated in Claim \ref{claim6}. Proceed in this way $b$ times. Finally, we obtain stars $S_1,\ldots,S_{b}$ with central vertices $z_1,\ldots,z_b.$ Let $R = [n]\setminus \left(\{1\}\cup N\cup V(S_1)\cup\ldots\cup V(S_b)\right).$ Due to the definition of $b$ and $\ell,$ $|R|<\floor[\big]{\frac{n}{\ln^2 n}}+m-2 < |N''|.$ Choose a subset $R'\subset N''$ of size $|R|.$ By Theorem \ref{bip}, whp there exists a matching $M$ between $R$ and $R'.$ Clearly, $\{S_1,\ldots,S_b\}\cup M\cup [N''\setminus R']$ is the desired $1$-factor.
\begin{flushright}
	$\blacksquare$
\end{flushright}

\noindent\textbf{Proof of Claim \ref{claim6}.} Let $\{u_{i,j}\}_{i,j=1}^{m-2,t}$ be vertices of $N$ (the labelling is arbitrary). Let $\tau_1\in S_t$ be identity. Let $B_1$ say that there exists a permutation $\tau_2\in S_t$ such that edges $\{u_{1,\tau_1(j)},u_{2,\tau_{2}(j)}\},$ $j\in[t],$ belong to $G(n,p).$ For every $i\in\{2,\ldots,m-4\},$ define $B_i$ recursively: $B_i$ says that there exists a permutation $\tau_{i+1}\in S_t$ such that edges $\{u_{i,\tau_i(j)},u_{i+1,\tau_{i+1}(j)}\},$ $j\in[t],$ belong to $G(n,p),$ and there are no edges between $u_{i+1,\tau_{i+1}(j)}$ and $u_{s,\tau_s(j)},$ $s\in[i-1],$ $j\in[t].$ Finally, let $B_{m-3}$ say that there exists a permutation $\tau_{m-2}\in S_t$ such that edges $\{u_{m-3,\tau_{m-3}(j)},u_{m-2,\tau_{m-2}(j)}\},$ $\{u_{1,\tau_{1}(j)},u_{m-2,\tau_{m-2}(j)}\},$ $j\in[t],$ belong to $G(n,p),$ and there are no edges between $u_{m-2,\tau_{m-2}(j)}$ and $u_{s,\tau_s(j)},$ $s\in\{2,\ldots,m-4\},$ $j\in[t].$

Consider a bipartite random graph $H_1$ with parts $U_1^1$ and $U_1^2,$ where vertices of $U_1^1$ are $u_{1,j},$ $j\in [t],$ vertices of $U_1^2$ are $u_{2,r},$ $r\in[t].$ Vertices are adjacent in this random bipartite graph iff they are adjacent in $G(n,p).$ Then $H_1$ is a binomial bipartite random graph with probability of appearing of an edge $p.$ Notice that $B_1$ holds iff $H_1$ has a perfect matching. The latter holds whp by Theorem \ref{bip}.

For every $i\in\{2,\ldots,m-4\},$ consider a random bipartite graph $H_i$ with parts $U_i^1$ and $U_i^2,$ where vertices of $U_i^1$ are sets $\{u_{1,\tau_{1}(j)},\ldots,u_{i,\tau_{i}(j)}\},$ $j\in[t],$ vertices of $U_i^2$ are $u_{i+1,r},$ $r\in[t].$ Vertices $\{u_{1,\tau_{1}(j)},\ldots,u_{i,\tau_{i}(j)}\}$ and $u_{i+1,r}$ are adjacent in this random bipartite graph iff
\begin{enumerate}
	\item $\{u_{i,\tau_i(j)},u_{i+1,r}\}\in E\left(G(n,p)\right),$
	\item there are no edges between $u_{i+1,r}$ and $u_{s,\tau_s(j)},$ $s\in[i-1],$ in $G(n,p).$
\end{enumerate}

Then $H_i$ is a binomial bipartite random graph with probability of appearing of an edge $p(1-p)^{i-1}.$ Notice that $B_i$ holds iff $H_i$ has a perfect matching. The latter holds whp by Theorem \ref{bip}.

Consider a bipartite random graph $H_{m-3}$ with parts $U^1_{m-3}$ and $U^2_{m-3},$ where vertices of $U_{m-3}^1$ are sets $\{u_{1,\tau_{1}(j)},\ldots,u_{m-3,\tau_{m-3}(j)}\},$ $j\in[t],$ vertices of $U_{m-3}^2$ are $u_{m-2,r},$ $r\in[t].$ Vertices $\{u_{1,\tau_{1}(j)},\ldots,u_{m-3,\tau_{m-3}(j)}\}$ and $u_{m-2,r}$ are adjacent in this random bipartite graph iff
\begin{enumerate}
	\item  $\{u_{1,\tau_1(j)},u_{m-2,r}\},$ $\{u_{m-3,\tau_{m-3}(j)},u_{m-2,r}\}\in E(G(n,p)),$
	\item there are no edges between $u_{m-2,r}$ and $u_{s,\tau_s(j)},$ $s\in\{2,\ldots,m-4\},$ in $G(n,p).$
\end{enumerate}
Then $H_{m-3}$ is a binomial bipartite random graph with probability of appearing of an edge $p^2(1-p)^{m-5}.$ Notice that $B_i$ holds iff $H_{m-3}$ has a perfect matching. The latter holds by Theorem \ref{bip}. Therefore, whp $G(n,p)[N]$ has $t$ vertex-disjoint cycles $C_{m-2}.$
\begin{flushright}
	$\blacksquare$
\end{flushright}

\noindent\textbf{Proof of Claim \ref{claim5}}. Let $X\subset N(1),$ $Y\subset[n]\setminus\left(\{1\}\sqcup N(1)\right)$ be sets of size $d:=\floor[\big]{\frac{n}{\ln^2 n}}.$ Let $\xi$ be the number of the desired $K_{1,a-1}.$ Then 
\begin{multline*}
	\lambda:={\sf E}\xi = d\binom{d}{a-1}p^{a-1}(1-p)^{\binom{a}{2}-a+1} = \\ = de^{(a-1)\ln\frac{d}{a-1}+a-a\ln\frac{1-p}{p}-\binom{a}{2}\ln\frac{1}{1-p}+O\left(\ln\ln n\right)}= e^{a\left(\ln\frac{d}{a}+1-\ln\frac{1-p}{p}-\frac{a-1}{2}\ln\frac{1}{1-p}\right)+O\left(\ln\ln n\right)}.
\end{multline*}
Since 
$$
\ln\frac{d}{a} - \frac{a}{2}\ln\frac{1}{1-p}=\ln n - 3\ln\ln n - \frac{2\log_{1/(1-p)}n - 8\log_{1/(1-p)}\ln n}{2}\ln\frac{1}{1-p}+O(1) = \ln\ln n + O(1),
$$ we get that $\lambda\geqslant e^{2(1+o(1))\ln\ln n\;\cdot\;\log_{1/(1-p)}n}.
$

Let $\xi = \xi_1+\ldots+\xi_{d\binom{d}{a-1}},$ where $\xi_i$ is the indicator of the presence of the $i$-th star in $G(n,p).$ Let $\Delta_{\ell} = \sum\limits_{i,j}{\sf E}\xi_i\xi_j$ over $i,j$ such that the respective stars have $\ell$ vertices in common. Then, for $\ell\geqslant 2,$
\begin{multline*}
	\Delta_{\ell}= d(d-1)\binom{d}{a-1}\binom{a-1}{\ell}\binom{d-a+1}{a-1-\ell}p^{2(a-1)}(1-p)^{2\binom{a}{2}-\binom{\ell}{2}-2a+2}+\\ d\binom{d}{a-1}\binom{a-1}{\ell-1}\binom{d-a+1}{a-\ell}p^{2(a-1)-(\ell-1)}(1-p)^{2\binom{a}{2}-\binom{\ell}{2}-2(a-1)+(\ell-1)}.
\end{multline*}
Let us put 
\begin{equation*}
	\Delta_{\ell}^1 := \frac{\binom{a-1}{\ell}\binom{d-a+1}{a-1-\ell}}{\binom{d}{a-1}}(1-p)^{-\binom{\ell}{2}},\quad \Delta_{\ell}^2 := \frac{\binom{a-1}{\ell-1}\binom{d-a+1}{a-\ell}}{d\binom{d}{a-1}}\left(\frac{1-p}{p}\right)^{\ell-1}(1-p)^{-\binom{\ell}{2}}.
\end{equation*}
It follows that we have
$$
\frac{\Delta_{\ell}}{\lambda^2}\leqslant \Delta_{\ell}^1+\Delta_{\ell}^2.
$$

Let us study the monotonicity of $\Delta_{\ell}^1$ and $\Delta_{\ell}^2$ on $[2, a-1].$ Notice that
\begin{equation*}
	\frac{\Delta_{\ell+1}^1}{\Delta_{\ell}^1} = \frac{(a-1-\ell)^2}{(\ell+1)(d-2a+\ell+3)}(1-p)^{-\ell},
\end{equation*}
 and
\begin{equation*}
	\frac{\partial}{\partial\ell}\ln\frac{\Delta_{\ell+1}^1}{\Delta_{\ell}^1} = -\frac{2}{a-1-\ell} - \frac{1}{\ell+1}-\frac{1}{d-2a+\ell+3}+\ln\frac{1}{1-p}.
\end{equation*}
Therefore, $\frac{\Delta_{\ell+1}^1}{\Delta_{\ell}^1}$ decreases when $\ell<\ell_1$ and $\ell>\ell_2$ and increases between $\ell_1,\ell_2,$ where $\ell_1 = O(1),$ $\ell_2=\ell-O(1)$ ($\ell_1$ may be equal to the left boundary of $[2,a-2]$).

Moreover, $\frac{\Delta_{\ell+1}^1}{\Delta_{\ell}^1}<1$ when $\ell=O(1)$ and $\frac{\Delta_{a-1}^1}{\Delta_{a-2}^1}>1.$ Therefore, there is a unique $\ell_0^1$ on $[2,a-2]$ such that $\frac{\Delta_{\ell_0^1+1}^1}{\Delta_{\ell_0^1}^1} = 1.$ We get that $\Delta_{\ell}$ decreases before $\ell_0^1$ and increases after $\ell_0^1.$

So, 
$$
\Delta_{\ell}^1\leqslant\max\{\Delta_2^1,\Delta_{a-1}^1\}=\max\left\lbrace \frac{\binom{a-1}{2}\binom{d-a+1}{a-3}}{\binom{d}{a-1}(1-p)},\frac{(1-p)^{-\binom{a-1}{2}}}{\binom{d}{a-1}} \right\rbrace.
$$
Observe that
\begin{equation*}
	\Delta_2^1\leqslant\frac{a^2\binom{d}{a-3}}{\binom{d}{a-1}(1-p)}<\frac{a^4}{d^2(1-p)}(1+o(1)),\quad \text{ and }\quad\Delta_{a-1}^1 = \frac{dp^{a-1}}{{\sf E}\xi}=e^{-2(1+o(1))\ln\ln n\cdot\log_{1/(1-p)}n}.
\end{equation*} 

In the same way, 
$$
\frac{\Delta_{\ell+1}^2}{\Delta_{\ell}^2} = \frac{(a-\ell)^2}{\ell(d-2a+\ell+2)}\left(\frac{1-p}{p}\right)(1-p)^{-\ell}
$$ 
equals to $1$ in a unique point $\ell_0^2,$ less than $1$ when $\ell<\ell_0^2$ and greater than $1$ when $\ell>\ell_0^2.$ So, for $n$ large enough,
$$
\Delta_{\ell}^2\leqslant\max\{\Delta_2^2,\Delta_{a-1}^2\}=\max\left\lbrace\frac{(a-1)\binom{d-a+1}{a-2}}{pd\binom{d}{a-1}},\frac{(a-1)(d-a+1)}{d\binom{d}{a-1}p^{a-2}(1-p)^{\binom{a-1}{2}-a+2}}\right\rbrace<\frac{a^2(1+o(1))}{d^2p}.
$$
By Janson's inequality \cite[Theorem 2.18]{JLR}, we obtain
\begin{multline*}
	{\sf P}(X=0)\leqslant \exp\left[-\frac{\lambda^2}{\sum\limits_{\ell=2}^a\Delta_{\ell}}\right]\leqslant \exp\left[-\frac{\lambda^2}{\lambda+a\max\{\Delta_2,\Delta_{\ell-1}\}}\right]
	=\\ \exp\left[-\frac{1}{\frac{1}{\lambda}+\frac{a^5(1+o(1))}{d^2(1-p)}}\right] = \exp\left[-\frac{d^2(1-p)(1+o(1))}{a^5}\right].
\end{multline*}

It remains to apply the union bound. By Theorem \ref{Chernoff}, whp $|N(1)|\in I,$ where 
$$
I=\left[ \floor[\big]{np-\sqrt{2n\ln n}},\floor[\big]{np+\sqrt{2n\ln n}}\right],
$$
we get that the probability of the complement to the desired event is at most 
$$
\max\limits_{s\in I}\binom{s}{d}\binom{n-s-1}{d}e^{-\Omega\left(\frac{n^2}{\ln^9n}\right)}+o(1)\to 0.
$$
\begin{flushright}
	$\blacksquare$
\end{flushright}
\section{Proof of Theorem \ref{th_main_2}}\label{proof_th3}
\begin{figure}
	\centering
	\includegraphics[width=0.6\linewidth]{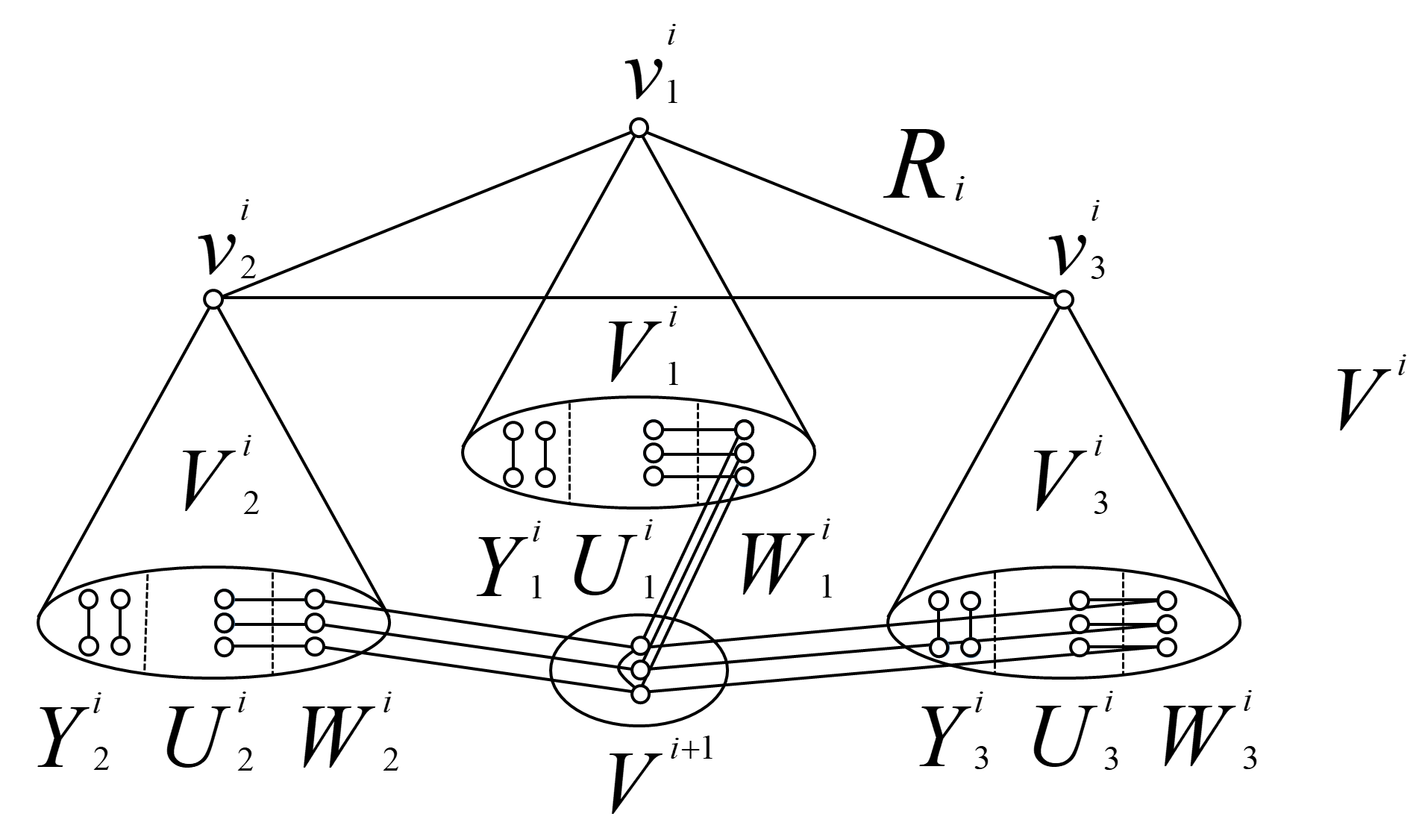}
	\caption{Graph $A[V^i]$}
	\label{fig:graphc4}
\end{figure}
To prove Theorem \ref{th_main_2}, we construct a $C_4$-saturated subgraph in $G(n,p)$ with the desired number of edges. The constructions for $p>1-1/\sqrt[3]{7}$ and $p\leqslant 1-\sqrt[3]{7}$ are different, we define them in Sections 4.1 and 4.2 respectively.
\subsection{$p>1-1/\sqrt[3]{7}$}
\begin{figure}
	\centering
	\includegraphics[width=0.6\linewidth]{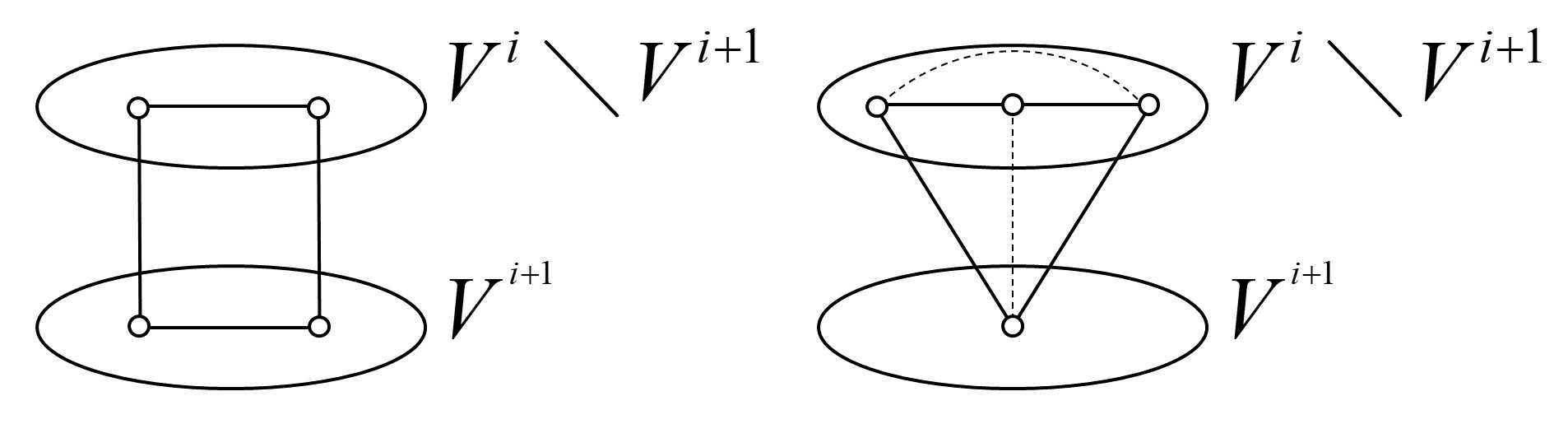}
	\caption{Dashed edges may be present.}
	\label{fig:graphconfig}
\end{figure}

Let $A$ be a graph on the vertex set $[n],$ $r\in\mathbb{N}.$ We call $A$ an $r$-flower (see Figure \ref{fig:graphc4}), if there exist sets of vertices $V^{r+1},V^{i},R_i,V_j^i,$ $j\in\{1,2,3\},$ $i\in[r],$ such that
\begin{enumerate}
	\item $V^{r+1}\sqcup\bigsqcup\limits_{i=\ell}^{r}\left(V^i_1\sqcup V^i_2 \sqcup V^i_3\sqcup R_i\right) = V^{\ell},$ $\ell \in[r],$\quad $V^1=[n];$
	\item for every $i \in[r],$ $R_{i}$ consists of pairwise adjacent (in A) vertices $v_1^{i},v_2^{i},v_3^{i},$ i.e. $A[R_{i}]$ is a copy of $K_3;$
	\item for every $j\in\{1,2,3\},$ $i\in[r],$  $N_{A[V^i]}(v_j^{i})\setminus R_i=V_j^i;$
	\item there are partitions $V_j^i=W_j^i\sqcup U_j^i\sqcup Y_j^i,$ $j\in\{1,2,3\},$ $i\in[r],$ such that
	\begin{enumerate}[label*=\text{\arabic*}]
		\item for every $i\in[r],$ $j\in\{1,2,3\},$ $E\left(A[W^i_j\sqcup V^{i+1}]\right)\setminus E\left(A[V^{i+1}]\right)$ is a perfect matching between $W^i_j$ and $V^{i+1},$
		\item for every $i\in[r],$ $j\in\{1,2,3\},$ $E\left(A[W^i_j\sqcup U_j^i] \right)$ is a perfect matching between $W^i_j$ and $U_j^i,$
		\item for every $i\in[r],$ $j\in\{1,2,3\},$ $E(A[Y_j^i])$ is a perfect matching;
	\end{enumerate}
	\item $A[V^{r+1}]$ is an inclusion-maximal $C_4$-free graph;
	\item there are no other edges in $A.$
\end{enumerate}
\begin{lemma}\label{lemma_Asat2}
	Let a graph $\Gamma$ contain a spanning $r$-flower $A.$ Then $A$ is $C_4$-saturated in $\Gamma$ and 
	\begin{equation}\label{edges_in_A}
		|E(A)|=|E(A[V^{r+1}])|+\frac{3}{2}\left(n-|V^{r+1}|-r\right)+3\left(|V^2|+\ldots+|V^{r+1}|\right).
	\end{equation}
\end{lemma}
\noindent\textbf{Proof.} First, we show that $A$ is $C_4$-free. 
\begin{itemize}
	\item For every $i\in[r],$ there are no copies of $C_4$ in $A[V^{i}\setminus V^{i+1}]$ since $A[V^{i}\setminus V^{i+1}]$ is a union of triangles that can be ordered in a way such that each triangle has a unique common vertex with the union of all its predecessors.
	\item $A[V^{r+1}]$ is $C_4$-free by the construction.
	\item Let $A[\{x_1,x_2,x_3,x_4\}]$ be a copy of $C_4.$ Let $i$ be the smallest number such that $V^i\setminus V^{i+1}$ contains one of the vertices $x_1,x_2,x_3,x_4.$ Every vertex $v$ of $V^i\setminus V^{i+1}$ is adjacent to at most one vertex in $V^{i+1}$ by the condition 4.1 of the definition of an $r$-flower. Therefore only two configurations remain (see Figure \ref{fig:graphconfig}). None of them is possible by the definition of an $r$-flower (conditions 2, 3, 4, 6).
\end{itemize}

Second, we show the maximality of $A.$ Let $i\in[r].$
\begin{itemize}
	\item For $j_1\neq j_2\in\{1,2,3\},$ an edge between $v_{j_1}^i$ and any $u\in V_{j_2}^i$ creates $C_4$ since $u$ has a neighbor $w\in V_{j_2}^i$ which is connected by an edge to $v_{j_2}^i.$ In turn, $v_{j_2}^i$ is connected to $v_{j_1}^i.$
	\item For any $j\in\{1,2,3\},$ an edge between $v_{j}^i$ and any $u\in V^{i+1}$ creates $C_4$ since there exists $w\in W_{j}^{i}$ connected with $u.$ In turn, $w$ has a neighbor in $U_j^i$ connected to $v_j^i.$
	\item For any $j_1\neq j_2\in\{1,2,3\},$ an edge between $u\in V^i_{j_1}$ and $v\in V^i_{j_2}$ creates $C_4$ since $u$ is connected to $v^i_{j_1},$ $v$ is connected to $v^i_{j_2}.$ In turn, $v^i_{j_1}$ and $v^i_{j_2}$ are adjacent.
	\item For any $j\in\{1,2,3\},$ an edge between any $u\in V_j^i$ and any $v\in V^{i+1}$ creates $C_4$ since there exists $w\in W_j^i$ connected to $v.$ In turn, both $w$ and $v$ are adjacent to $v_j^i.$
\end{itemize}

Finally, $A[V^{r+1}]$ is maximal by the definition.\\

It stays on to verify \eqref{edges_in_A}. Removal of $v_1^i$ and of edges adjacent to $v_1^i,$ $v_2^i,$ $v_3^i$ except for $\{v_2^i,v_3^i\},$ makes $E\left(A[V^i\setminus V^{i+1}]\right)$ a matching for every $i\in[r].$ The number of the removed edges is $n-r-|V^{r+1}|$ while the cardinality of the matching is half this number. It is only left to take into account $3|V^{i+1}|$ edges between $V^i$ and $V^{i+1}$ as well as $|E(A[V^{r+1}])|$ edges in $A[V^{r+1}].$
\begin{flushright}
	$\blacksquare$
\end{flushright}

\begin{lemma}\label{lemma_property_t}
	For every $p>1-1/\sqrt[3]{7},$ whp $G(n,p)$ contains a spanning $r$-flower $A$ where $r = \ceil[\big]{\frac{5}{24}\log_{\frac{1}{1-p}}n},$ $|V^2|+\ldots+|V^{r+1}| = \frac{n((1-p)^3+o(1))}{1-(1-p)^3}$ and $|V^{r+1}|=o\left(\sqrt{n}\right).$
\end{lemma}
\noindent\textbf{Proof.} Below we recursively construct sets $[n]=V^1\reflectbox{ $\subset$ }\ldots V^{r}\reflectbox{ $\subset$ }V^{r+1}.$ We use the following notations for some characteristics of graphs induced by these sets in $G(n,p).$ For any $s\in\mathbb{N},$ $i\in[r],$ $v_1,\ldots,v_s\in V^i$ and $\nu_k\in\{v_k,\overline{v}_k\},$ $k\in[s],$ let us denote by $N_i(\nu_1,\ldots,\nu_{s})$ the set of all common neighbors of $\{v_k,\; k\in[s]\;|\;\nu_k = v_k\}$ among common non-neighbors of $\{v_k,\; k\in[s]\;|\;\nu_k = \overline{v}_k\}$ in $V^i$ (not including $v_k \text{ such that }\nu_k = \overline{v}_k$). Let $n_i(\nu_1,\ldots,\nu_{s})$ stand for the cardinality of $N_i(\nu_1,\ldots,\nu_{s}).$

Put $\varepsilon=(1/(1-p)^3-7)/4$ and set $V^1=[n].$ Notice that $\varepsilon>0$ due to the restriction on $p.$ Assume that, for some $i\in[r],$ the sets $V^1\reflectbox{ $\subset$ }\ldots \reflectbox{ $\subset$ }V^{i}$ are already constructed and that $G(n,p)[V^i]\stackrel{d}=G(|V^i|,p).$ Suppose that with probability $1-O\left(i/\sqrt{n}\right)$ (the constant factor in $O(i/\sqrt{n})$ does not depend on $i$), 
\begin{equation} \label{Vi}
	|V^j|=  n(1-p)^{3(j-1)}+O\left(j\sqrt{n(1-p)^{3(j-1)}\ln n}\right),\quad\text{for all } j\in[i].
\end{equation} 
Since $i\leqslant r,$ this implies $|V^i|\gg n^{1/3}.$ All the probabilities below are conditioned on the event \eqref{Vi}. Choose an arbitrary vertex $v_1^i$ in $V^i.$ Denote by $B$ the event that 
\begin{itemize}
	\item $|V^i|p-\sqrt{|V^i|p\ln n}\leqslant n_i(v_1^i)\leqslant |V^i|p+\sqrt{|V^i|p\ln n},$ 
	\item $|V^i|(1-p)-\sqrt{|V^i|(1-p)\ln n}\leqslant n_i(\overline{v}_1^i)\leqslant|V^i|(1-p)+\sqrt{|V^i|(1-p)\ln n}.$
\end{itemize}
By Theorem \ref{Chernoff}, the probability that $B$ holds is $1-O\left(1/\sqrt{n}\right).$ $B$ (jointly with \eqref{Vi}) implies $n_i(v_1^i)\gg n^{1/3}$ and $n_i(\overline{v}_1^i)\gg n^{1/3}.$ Expose edges adjacent to $v_1^i$ in $V^i.$ Choose any neighbor of $v_1^i$ (if exists) and denote it by $v_2^i.$ Denote by $D$ the event that
\begin{itemize}
	\item $n_i(v_1^i)p-\sqrt{n_i(v_1^i)p\ln n} \leqslant n_i(v_1^i,v_2^i)\leqslant n_i(v_1^i)p-\sqrt{n_i(v_1^i)p\ln n},$
	\item $n_i(v_1^i)(1-p)-\sqrt{n_i(v_1^i)(1-p)\ln n}\leqslant n_i(v_1^i,\overline{v}_2^i)\leqslant n_i(v_1^i)(1-p)+\sqrt{n_i(v_1^i)(1-p)\ln n},$
	\item $n_i(\overline{v}_1^i)p-\sqrt{n_i(\overline{v}_1^i)p\ln n}\leqslant n_i(\overline{v}_1^i,v_2^i)\leqslant n_i(\overline{v}_1^i)p+\sqrt{n_i(\overline{v}_1^i)p\ln n},$
	\item $n_i(\overline{v}_1^i)(1-p)-\sqrt{n_i(\overline{v}_1^i)(1-p)\ln n}\leqslant n_i(\overline{v}_1^i,\overline{v}_2^i)\leqslant n_i(\overline{v}_1^i)(1-p)+\sqrt{n_i(\overline{v}_1^i)(1-p)\ln n}.$
\end{itemize}
By Theorem \ref{Chernoff}, the probability that $D\cap B$ holds is $1-O\left(1/\sqrt{n}\right).$ These inequalities (jointly with \eqref{Vi} and $B$) imply that $n_i(v_1^i,v_2^i)\gg n^{1/3},$ $n_i(v_1^i,\overline{v}_2^i)\gg n^{1/3},$ $n_i(\overline{v}_1^i,v_2^i)\gg n^{1/3}$ and $n_i(\overline{v}_1^i,\overline{v}_2^i)\gg n^{1/3}.$ Expose edges adjacent to $v_2^i$ in $V^i.$ Choose an arbitrary common neighbor of $v_1^i$ and $v_2^i$ (if exists) and denote it by $v_3^i.$ In the same way, denote by $K$ the event that each $n_i(\nu_1,\nu_2,\nu_3),$ $\nu_j\in\{v_j^i,\overline{v}_j^i\},$ differs from its expectation ${\sf E} n_i:={\sf E} \left(n_i(\nu_1,\nu_2,\nu_3)|n_i(\nu_1,\nu_2)\right)$ by at most $\sqrt{{\sf E}n_i\ln n}.$ By Theorem \ref{Chernoff}, the probability that $K\cap D\cap B$ holds is $1-O\left(1/\sqrt{n}\right).$ Expose edges adjacent to $v_3^i$ in $V^i.$\\

Choose arbitrarily
\begin{itemize}
	\item $\ceil{n_i(v_1^i,v_2^i,v_3^i)/3}$ vertices of $N_i(v_1^i,v_2^i,v_3^i),$
	\item $\ceil{n_i(v_1^i,\overline{v}_2^i,v_3^i)/2}$ vertices of $N_i(v_1^i,\overline{v}_2^i,v_3^i),$
	\item $\ceil{n_i(v_1^i,v_2^i,\overline{v}_3^i)/2}$ vertices of $N_i(v_1^i,v_2^i,\overline{v}_3^i),$
	\item all vertices of $N_i(v_1^i,\overline{v}_2^i,\overline{v}_3^i)$ 
\end{itemize}
into $V^i_1.$ Choose arbitrarily 
\begin{itemize}
	\item another $\ceil{n_i(v_1^i,v_2^i,v_3^i)/3}$ vertices of $N_i(v_1^i,v_2^i,v_3^i),$
	\item the remaining vertices of $N_i(v_1^i,v_2^i,\overline{v}_3^i),$
	\item $\ceil{n_i(\overline{v}_1^i,v_2^i,v_3^i)/2}$ vertices of $N_i(\overline{v}_1^i,v_2^i,v_3^i),$
	\item all vertices of $N_i(\overline{v}_1^i,v_2^i,\overline{v}_3^i)$
\end{itemize}
into $V^i_2.$ Choose the remaining vertices from $N_i(v_1^i,v_2^i,v_3^i),$ $N_i(\overline{v}_1^i,v_2^i,v_3^i),$ $N_i(v_1^i,\overline{v}_2^i,v_3^i)$ and all the vertices of $N_i(\overline{v}_1^i,\overline{v}_2^i,v_3^i)$ into $V^i_3.$ It follows from $K\cap D\cap B$ that 
$$
|V_j^i|=|V^i|(1-(1-p)^3)/3+O\left(\sqrt{|V^i|\ln n}\right)\quad\text{ for every } j\in\{1,2,3\}
$$
and 
\begin{equation}\label{ni_up}
	n_i(\overline{v}_1^i,\overline{v}_2^i,\overline{v}_3^i)= |V^i|(1-p)^3+O\left(\sqrt{|V^i|\ln n}\right).
\end{equation}
It implies that
\begin{equation}\label{eq_eps}
	|V_j^i|-2n_i(\overline{v}_1^i,\overline{v}_2^i,\overline{v}_3^i)>\varepsilon n_i(\overline{v}_1^i,\overline{v}_2^i,\overline{v}_3^i).
\end{equation}

For every $j\in\{1,2,3\},$ choose arbitrarily $|V_j^i|-2n_i(\overline{v}_1^i,\overline{v}_2^i,\overline{v}_3^i)$ vertices from $V_j^i.$ Denote the sets of the chosen vertices by $A_j.$ In every $A_j,$ choose $2\tau$ vertices and add them to $V_j^i\setminus A_j,$ where $\tau$ stands for the number of $A_j$ with odd cardinalities. In every $A_j$ with an odd cardinality choose any vertex, remove edges joining it to $v_j^i$ and add this vertex to $N_i(\overline{v}_1^i,\overline{v}_2^i,\overline{v}_3^i).$ Split arbitrarily the sets obtained from $V_j^i\setminus A_j$ into two equal parts and denote the parts by $W_j^i$ and $U_j^i.$ Designate by $Y_j^i$ the set $V_j^i\setminus\left(U_j^i\sqcup W_j^i\right).$ We denote by $V^{i+1}$ the obtained set $N_i(\overline{v}_1^i,\overline{v}_2^i,\overline{v}_3^i)$ (with added $\tau$ vertices). Therefore,
\begin{equation}\label{Vi_ni}
	|V^{i+1}| = n_i(\overline{v}_1^i,\overline{v}_2^i,\overline{v}_3^i) + O(1)
\end{equation}
(the $O(1)$ is bounded from above uniformly over all $i$). Expression \eqref{ni_up} (jointly with \eqref{Vi}) implies $n_i(\overline{v}_1^i,\overline{v}_2^i,\overline{v}_3^i)\gg n^{1/3}.$

By \eqref{er}, the probability that $G(n,p)[Y_j^i]$ has no perfect matching, given $n_i(\overline{v}_1^i,\overline{v}_2^i,\overline{v}_3^i)\gg n^{1/3}$ and \eqref{eq_eps}, is $e^{-\omega(n^{1/3})}.$ Observe that the graph induced by $G(n,p)$ between $W_j^i$ and $V^{i+1}$ is a bipartite random graph with the probability of an edge $p.$ Similarly, the graph induced by $G(n,p)$ between $W_j^i$ and $U_j^i$ is a bipartite random graph with the probability of an edge $p.$ By Theorem \ref{bip}, the probability that either there is no perfect matching between $W_j^i$ and $V^{i+1}$ or there is no perfect matching between $W_j^i$ and $U_j^i,$ given  $n_i(\overline{v}_1^i,\overline{v}_2^i,\overline{v}_3^i)\gg n^{1/3},$ is $e^{-\omega(n^{1/3})}.$ Expose edges adjacent to $W_j^i,$ $Y_j^i,$ $U_j^i$ in $V^i.$ Find perfect matchings between $W_j^i$ and $V^{i+1},$ between $W_j^i$ and $U_j^i,$ in $Y_j^i$ (if exist).

Remove all edges in $V^i\setminus V^{i+1}$ except for those between $v_1^i,$ $v_2^i$ and $v_3^i;$ between $v^i_j$ and $V^i_j;$ in the perfect matchings found in $Y_j^i;$ in the perfect matchings found between $U_j^i$ and $W_j^i.$ Remove all edges between $V^i\setminus V^{i+1}$ and $V^{i+1}$ except for the perfect matchings found between $W_j^i$ and $V^{i+1}.$ Evidently, $G(n,p)[V^{i+1}]\stackrel{d}=G(|V^{i+1}|,p).$ Finally, from \eqref{Vi} -- \eqref{Vi_ni}, we get that with probability $1-O((i+1)/\sqrt{n}),$
\begin{equation}\label{eqip1}
	|V^{i+1}| =  n(1-p)^{3i}+O\left((i+1)\sqrt{n(1-p)^{3i}\ln n}\right)
\end{equation}
(where the big-O is bounded by the same constant as in \eqref{Vi}). This finishes the step of the induction. Since $r\ll \sqrt{n},$ whp we construct all $V^1\reflectbox{ $\subset$ }\ldots \reflectbox{ $\subset$ }V^{r+1}$ successfully.\\

Expose edges within $V^{r+1}.$ Remove some edges such that a $C_4$-saturated graph is left.

It remains to estimate $|V^2|+\ldots+|V^r|$ and $|V^{r+1}|.$ It follows from \eqref{eqip1} that whp
$$
|V^{r+1}|=n(1-p)^{3r}(1+o(1))=o(\sqrt{n}),
$$
\begin{equation*}
	|V^2|+\ldots+|V^{r+1}| = \frac{n(1-p)^3}{1-(1-p)^3}+o(n).
\end{equation*}
\begin{flushright}
	$\blacksquare$
\end{flushright}

Due to Lemma \ref{lemma_Asat2} and Lemma \ref{lemma_property_t}, for every $\varepsilon>0,$ whp $G(n,p)$ contains a subgraph $A$ with at most 
$$
\frac{3(1+(1-p)^3)}{2(1-(1-p)^3)}n(1+o(1))
$$
edges such that $A$ is $C_4$-saturated. Hence, whp
$$
\mathrm{sat}\left(G(n,p),C_4\right)\leqslant \frac{3(1+(1-p)^3)}{2(1-(1-p)^3)}n(1+o(1)).
$$

The inequality \eqref{cycle_4_up} is proven.

\subsection{$p\leqslant 1-1/\sqrt[3]{7}$}
\begin{figure}
	\centering
	\includegraphics[width=0.8\linewidth]{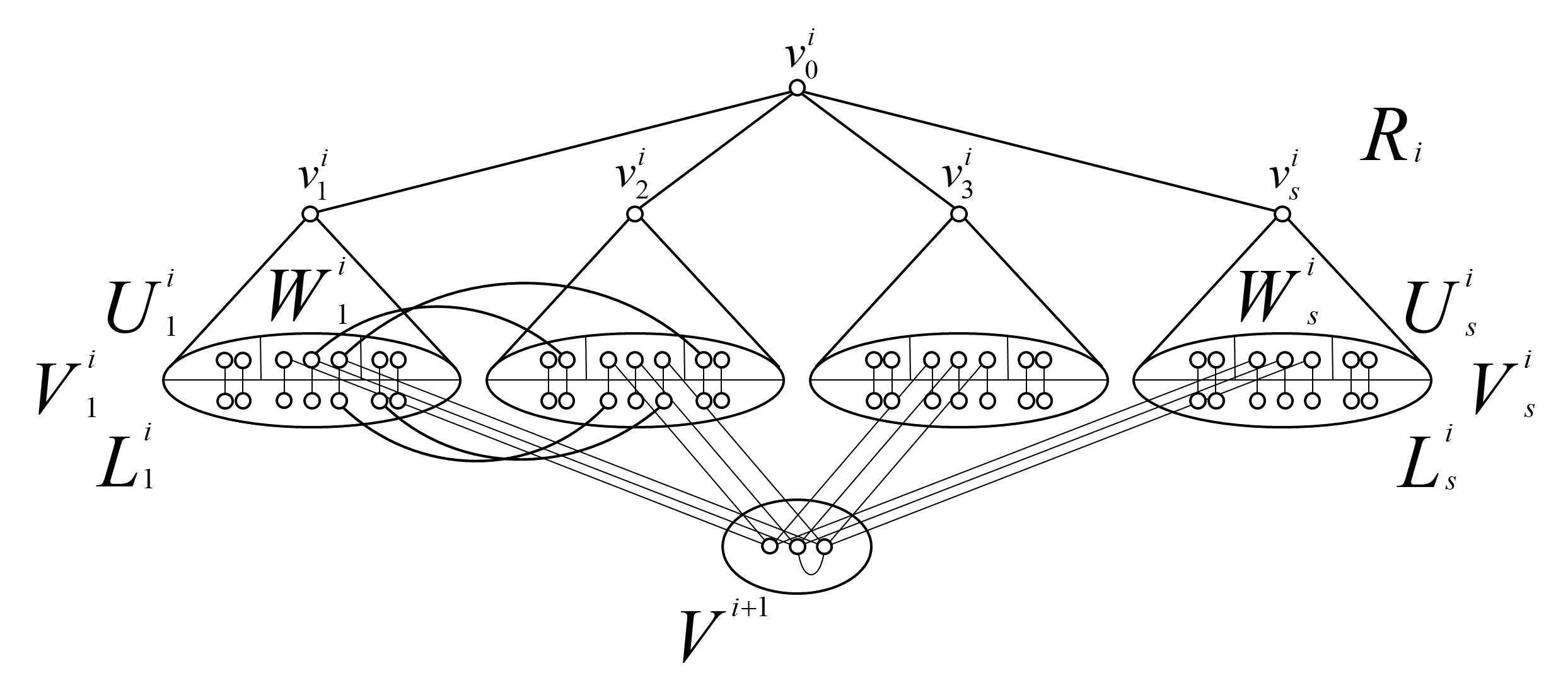}
	\caption{Graph $A[V^i].$}
	\label{fig:graph_prop_M}
\end{figure}

\begin{figure}
	\centering
	\begin{subfigure}{0.33\textwidth}
		\centering
		\includegraphics[height=3.5cm]{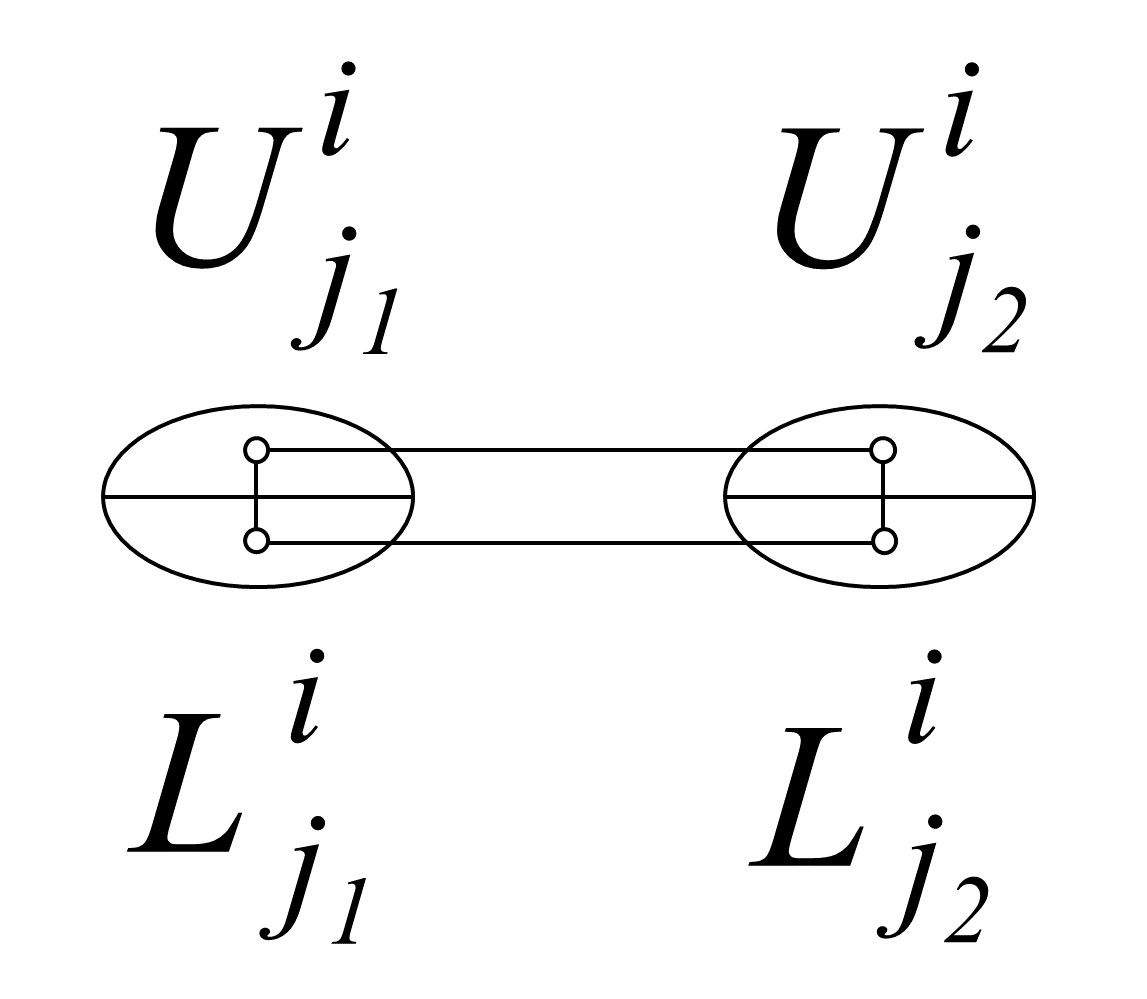}
		\caption{}
		\label{fig:3}
	\end{subfigure}
	\begin{subfigure}{0.33\textwidth}
		\centering
		\includegraphics[height=3.5cm]{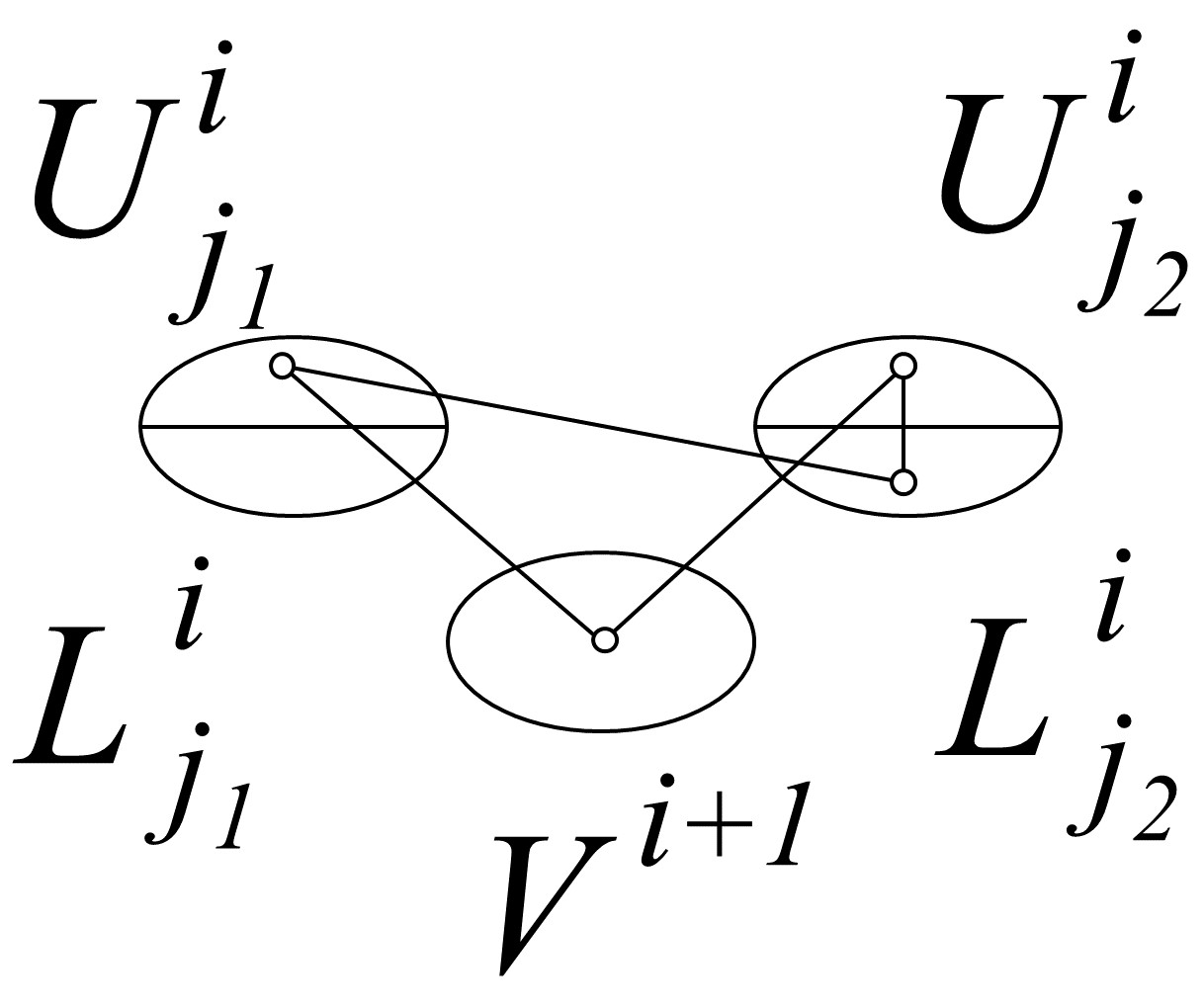}
		\caption{}
		\label{fig:1}
	\end{subfigure}
	\begin{subfigure}{0.33\textwidth}
		\centering
		\includegraphics[height=3.5cm]{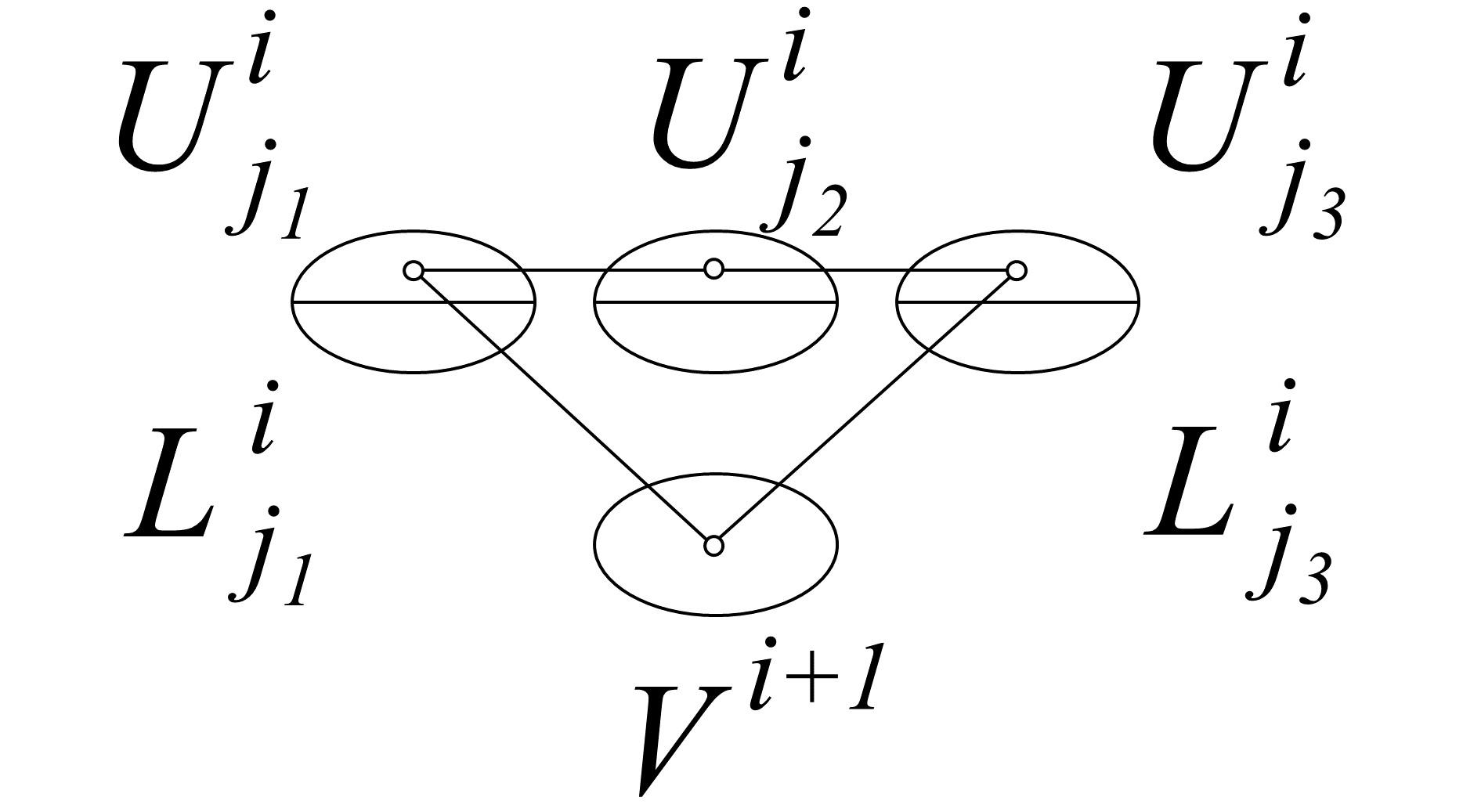}
		\caption{}
		\label{fig:2}
	\end{subfigure}
	\begin{subfigure}{0.33\textwidth}
		\centering
		\includegraphics[height=3.5cm]{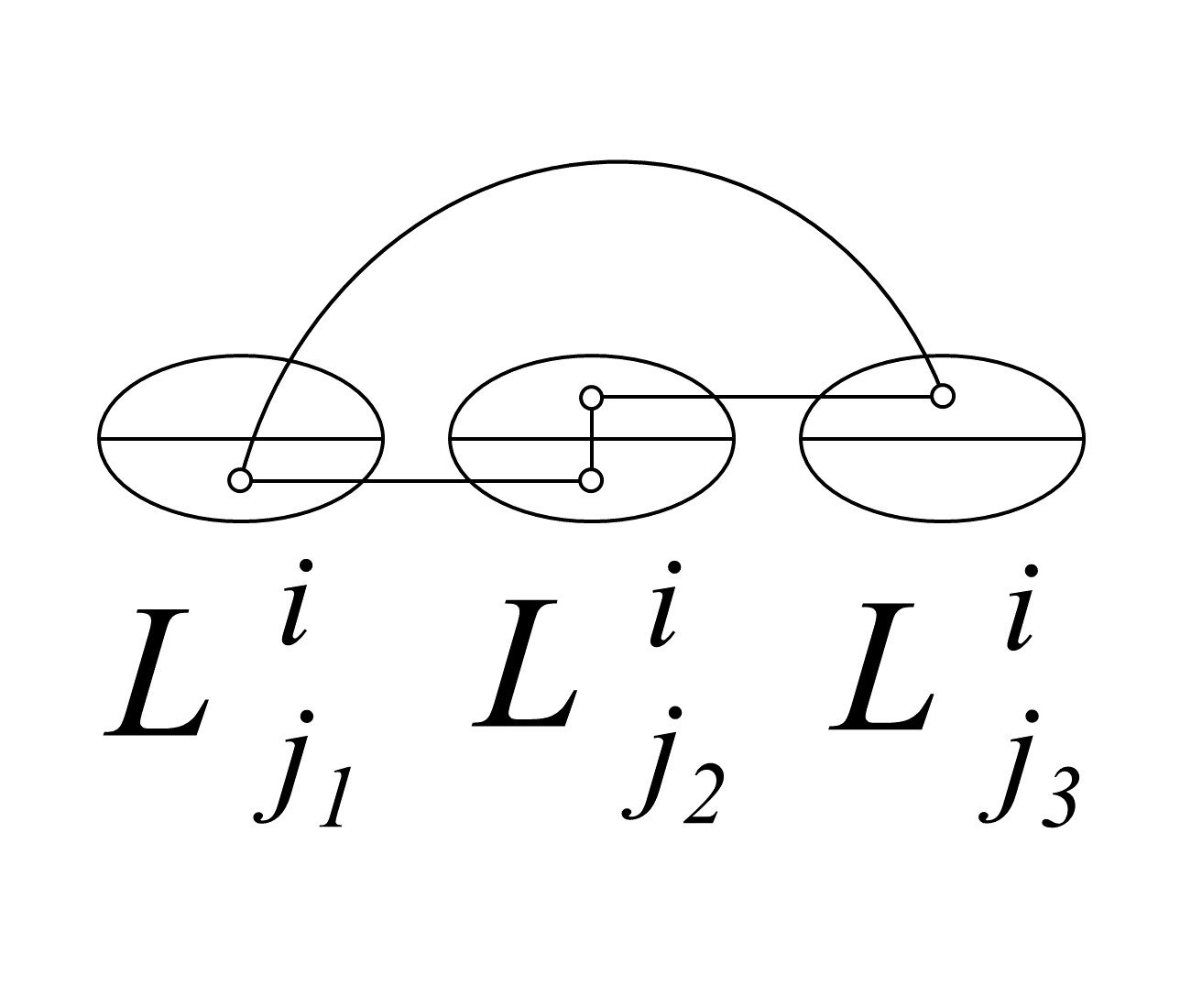}
		\caption{}
		\label{fig:4}
	\end{subfigure}\\
	\begin{subfigure}{0.33\textwidth}
		\centering
		\includegraphics[height=3cm]{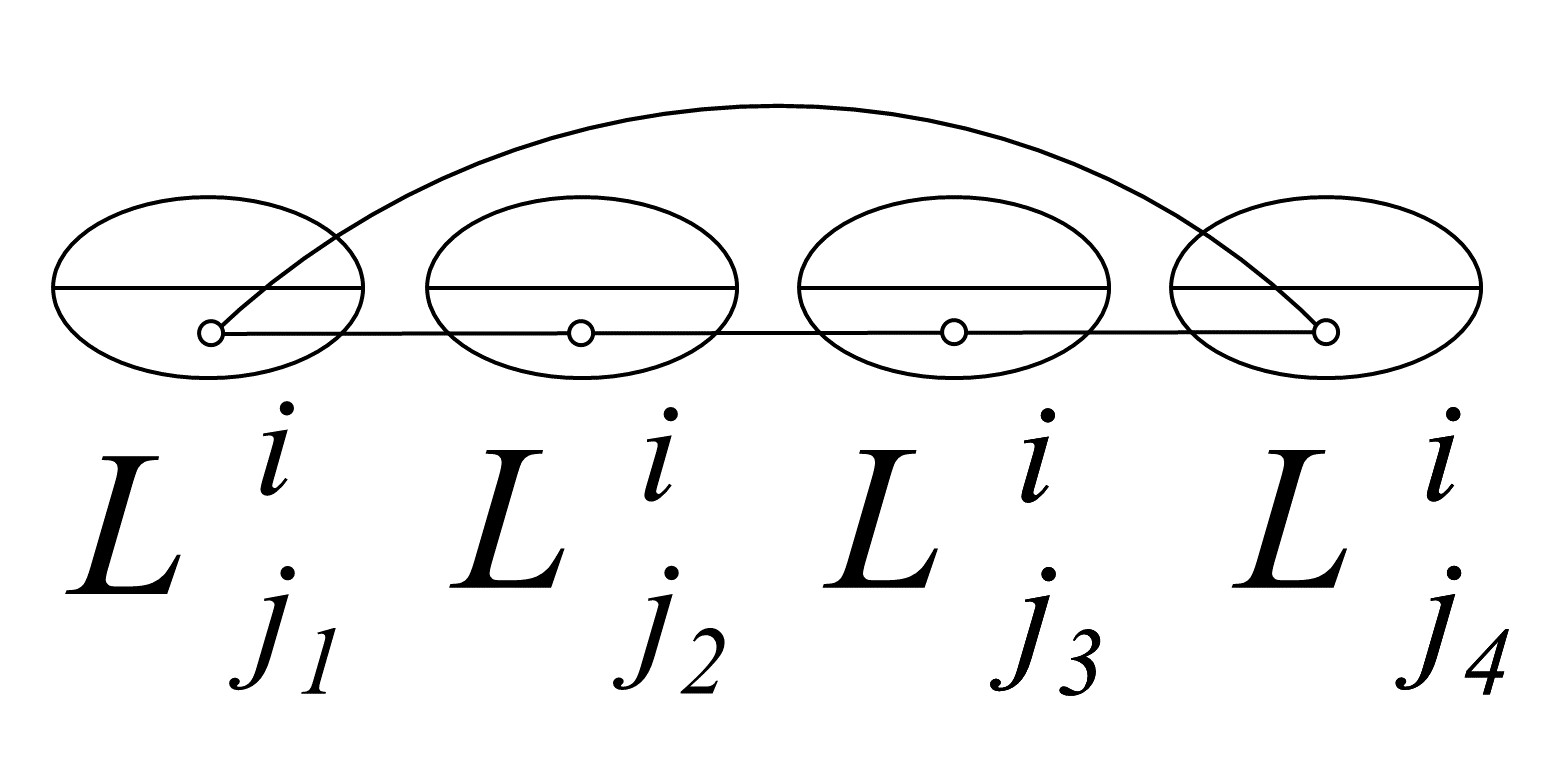}
		\caption{}
		\label{fig:5}
	\end{subfigure}
	\begin{subfigure}{0.33\textwidth}
		\centering
		\includegraphics[height=3cm]{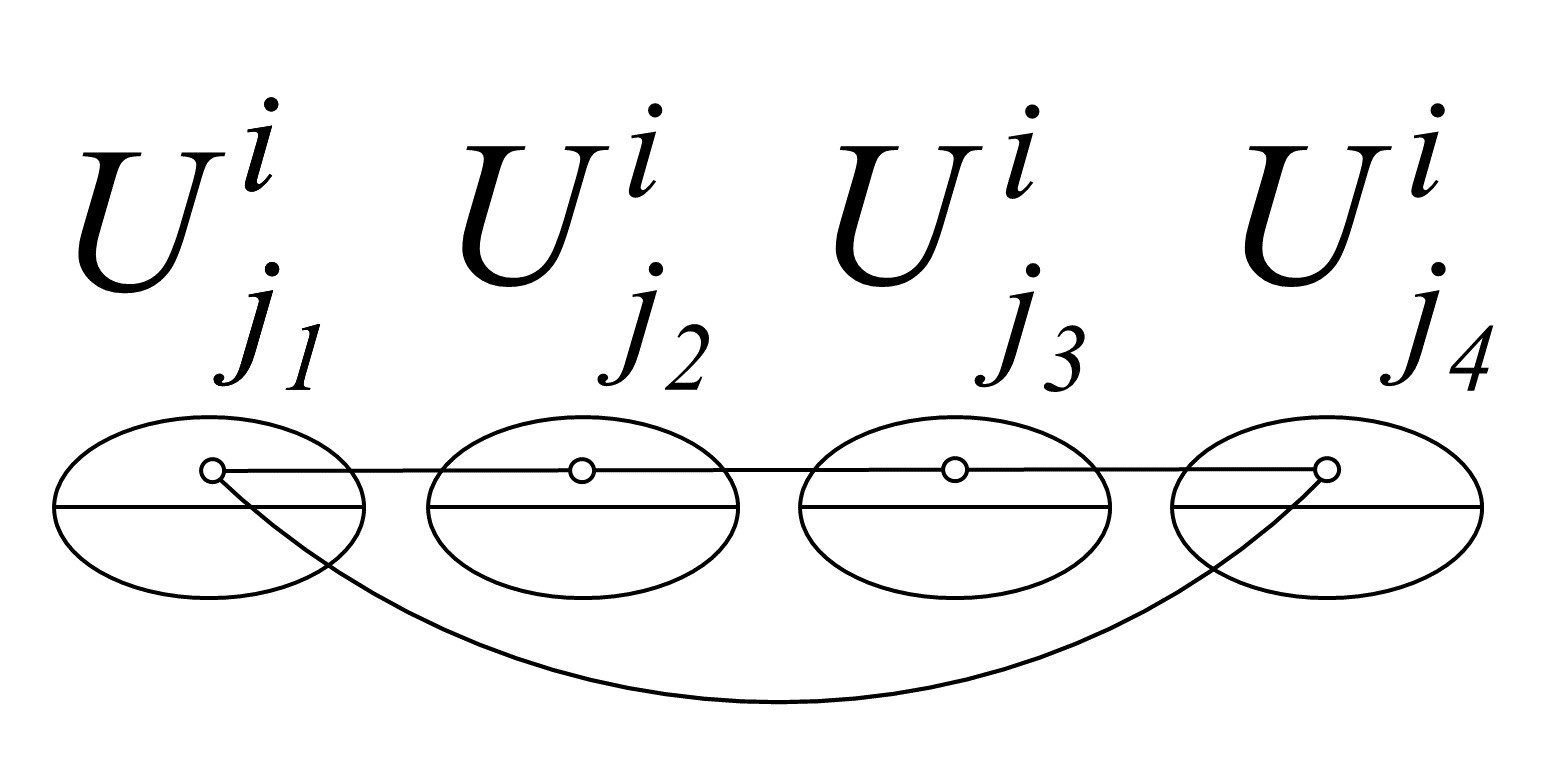}
		\caption{}
		\label{fig:7}
	\end{subfigure}\\
	\begin{subfigure}{0.33\textwidth}
		\centering
		\includegraphics[height=3.5cm]{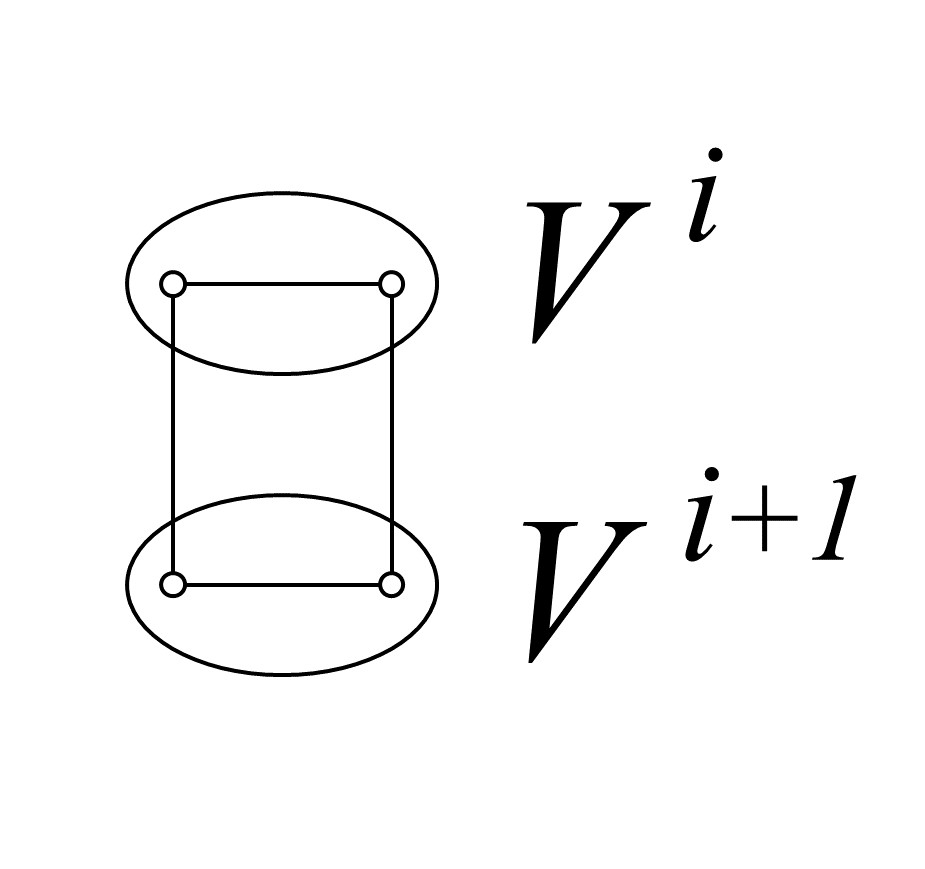}
		\caption{}
		\label{fig:6}
	\end{subfigure}
	\caption{}
\end{figure}

Let $A$ be a graph on the vertex set $[n],$ $r,s\in\mathbb{N}.$ We call $A$ an $(s,r)$\textit{-flower} (see Figure \ref{fig:graph_prop_M}), if there exist sets of vertices $V^{r+1},V^i,R_i,V_j^i,$ $j\in[s],$ $i\in[r],$ such that
\begin{enumerate}
	\item $V^{r+1}\sqcup \bigsqcup\limits_{i=\ell}^r\left(V_1^i\sqcup\ldots\sqcup V_s^i\sqcup R_i\right)=V^{\ell},\;\ell\in[r],\quad V^1=[n];$
	\item for every $i\in[r],$ $R_i=\{v_0^i,v_1^i,\ldots,v_s^i\};$
	\item for every $i\in[r],$ $j\in[s],$ $\{v_1^i,\ldots,v^i_s\}= N_{A[V^i]}(v^i_0),$ $V_j^i=N_{A[V^i]}(v^i_j)\setminus\{v_0^i\};$
	\item there are partitions $V^i_j=U_{j}^i\sqcup L_{j}^i$ and subsets $W_{j}^i\subset U_{j}^i,$ $i\in[r],$ $j\in[s],$ such that 
	\begin{enumerate}[label*=\text{\arabic*}]
		\item $|U_{j}^i|=|L_{j}^i|,$ $|W_j^i|=|V^{i+1}|;$
		\item $E\left(A\left[U_{j}^i\sqcup L_{j}^i\right]\right)$ is a perfect matching between $U_{j}^i$ and $L_{j}^i;$
		\item $E(A[W_{j}^i\sqcup V^{i+1}])\setminus E(A[V^{i+1}])$ is a perfect matching between $W_{j}^i$ and $V^{i+1};$
		\item for every $j_1\neq j_2\in[s],$   $E\left(A[U_{j_1}^{i}\sqcup U_{j_2}^{i}]\right)$ is a perfect matching between $U_{j_1}^{i}$ and $U_{j_2}^{i};$
		\item for every $j_1\neq j_2\in[s],$   $E\left(A[L_{j_1}^{i}\sqcup L_{j_2}^{i}]\right)$ is a perfect matching between $L_{j_1}^{i}$ and $L_{j_2}^{i};$
		\item for every $j_1\neq j_2\in[s],$ $E\left(A[W_{j_1}^{i}\sqcup W_{j_2}^{i}]\right)$ is empty;
		\item copies of $C_4$ presented on Figures \ref{fig:3}, \ref{fig:2}, \ref{fig:5}, \ref{fig:7} do not appear in $A;$
	\end{enumerate}
	\item $A[V^{r+1}]$ is an inclusion-maximal $C_4$-free graph;
	\item there are no other edges in $A.$
\end{enumerate}
\begin{lemma}\label{lemma_Asat3} Let $\Gamma$ contain a spanning $(s,r)$-flower $A.$ Then $A$ is $C_4$-saturated in $\Gamma$ and 
\begin{equation}\label{eq:edges_init}
	|E(A)|=|E(A[V^{r+1}])|-r+\frac{s+1}{2}\left(n-|V^{r+1}|\right)-\frac{rs(s+1)}{2}+s\left(|V^2|+\ldots+|V^{r+1}|\right).
\end{equation}
\end{lemma}
\noindent\textbf{Proof.} First, we show that $A$ is $C_4$-free. 

Let us show that for every $i\in[r],$ there are no copies of $C_4$ in $A[V^i\setminus V^{i+1}].$ The vertex $v_0^i$ can not be in any $C_4\subset A[V^i\setminus V^{i+1}]$ since none of its neighbors have other common neighbors by the definition of an $(s,r)$-flower (conditions 1, 2, 3, 6). It  implies that any vertex $v_j^i$ may only be contained in $C_4\subset A[V^i\setminus V^{i+1}]$ with both neigbors from $V_j^i.$ However, by the definition of an $(s,r)$-flower, no two vertices from $V_j^i$ have a common neighbor other than $v_j^i$ (conditions 3, 4.2, 4.4, 4.5, 6). Also, by the definition,
\begin{itemize}
	\item no vertex from $V_{j}^i,$ $j\in[s],$ can have two neighbors in $V^{i}\setminus(R_i\sqcup V_j^i\sqcup V^{i+1})$ having another common neighbor in $V^i\setminus V^{i+1}$ (conditions 4.4, 4.5, 4.6, 6);
	\item no two connected vertices from $V_j^i,$ $j\in[s],$ can be in $C_4\subset A[V^i\setminus V^{i+1}]$ (see Figures \ref{fig:3}, \ref{fig:4}) (conditions 4.2, 4.4, 4.5, 4.6, 6).
\end{itemize}
Also, there are no copies of $C_4$ having vertices both in $V^i$ and $V^{i+1}$. Indeed, no vertex from $V_j^i, $ $j\in[s],$ can have two neighbors in $V^{i+1}$ (conditions 4.3, 6). A configuration shown on Figure \ref{fig:2} is not possible by the condition 4.6. Since $U$-sets have only neighbors in $U$-sets, and $L$-sets have only neighbors in $L$-sets (conditions 4.4, 4.5, 6), a configuration given on Figure \ref{fig:1} is not possible. A configuration demonstrated on Figure \ref{fig:6} is not possible by the conditions 4.3, 4.6, 6. The vertex $v_0^i$ can not be in $C_4$ with vertices from $V^{i+1}$ since neither $v_0^i$ nor its neigbors $v^i_j,$ $j\in[s],$ are in edges adjacent to $V^{i+1}$ (conditions 3, 6). No $v_j^i$ is contained in a copy of $C_4$ with vertices from $V^{i+1}$ since its neighbors belong to $V^i\setminus V^{i+1},$ and do not have common neighbors in $V^{i+1}$ (conditions 3, 4.3, 6).

Finally, $A[V^{r+1}]$ is $C_4$-free by the condition 5.\\

Second, we prove the maximality of $A.$ Let $i\in[r].$
\begin{itemize}
	\item An edge between $v_0^i$ and any $v\in V^{i+1}$ creates $C_4$ since there is a vertex $u\in W_{1}^i$ connected to $v.$ In turn, both $u$ and $v_0^i$ are adjacent to $v_1^i.$
	\item An edge between $v_0^i$ and any $v\in V_{j}^i$ creates $C_4$ since $E\left(A[U_{j}^i\sqcup L_{j}^i]\right)$ is a perfect matching. Indeed, if $\{u,v\}$ is an edge of this perfect matching, then $u\sim v_j^i,$ $v_j^i\sim v_0^i.$
	\item For $j_1\neq j_2\in[s],$ an edge between $v_{j_1}^i$ and $v_{j_2}^i$ creates $C_4$ since $E\left(A[U_{j_1}^i\sqcup U_{j_2}^i]\right)$ is a perfect matching. Indeed, if $\{u,v\}$ is an edge of this perfect matching, then $u\sim v_{j_1}^i,$ $v\sim v_{j_2}^i.$
	\item For $j_1\neq j_2\in[s],$ an edge between $v_{j_1}^i$ and any $u\in V_{j_2}^i$ creates $C_4$ since $u$ is connected to $v_{j_2}^i.$ In turn, both $v_{j_1}^i$ and $v_{j_2}^i$ are adjacent to $v_0^i.$
	\item For $j\in[s],$ an edge between $v_j^i$ and any $w\in V^{i+1}$ creates $C_4$ since there exists $u\in W_{j}^i$ connected to $w.$ In turn, there exists $v\in L_{j}^i$ adjacent to both $u$ and $v_j^i.$
	\item For $j\in[s],$ an edge between $v\in V_j^i$ and any $w\in V^{i+1}$ creates $C_4$ since there exists $u\in W_{j}^i$ connected to $w.$ In turn, both $u$ and $v$ are adjacent to $v_j^i.$
	\item For $j\in [s],$ an edge between $u$ an $v$ from $V^i_j$ creates $C_4$ since there exists $w\in V^i_j$ such that $u\sim w$ and $w\sim v^i_j,$ $v\sim v^i_j.$
	\item For $j_1\neq j_2\in[s],$ an edge between $u\in V^i_{j_1}$ and $v\in V^i_{j_2}$ creates $C_4$ since there exists $w\in V^i_{j_1}$ such that $v\sim w.$ In turn, $u\sim v^i_{j_1}$ and $w\sim v^i_{j_1}.$
\end{itemize}
Finally, $A[V^{r+1}]$ is maximal by the condition 6.\\

It remains to count edges in $A.$ There are 
\begin{itemize}
	\item $|E(A[V^{r+1}])|$ edges in $A[V^{r+1}];$
	\item $n-1-|V^{r+1}|-(r-1)$ edges in spanning trees
	$$
	\left(V^i\setminus V^{i+1},E\left(A[R_i]\right)\sqcup\bigsqcup_{j=1}^s \left(E\left(A[V_j^i\cup v_j^i]\right)\setminus E\left(A[V_j^i]\right)\right)\right)
	$$
	of $A[V^i\setminus V^{i+1}],$ $i\in[r];$
	\item $\left(n-\sum_{i=1}^{r}|R_i|-|V^{r+1}|\right)/2$ edges in
	\begin{equation*}
		\bigsqcup\limits_{i=1}^r\bigsqcup_{j=1}^s E\left(A[V_j^i]\right);
	\end{equation*} 
	\item $\left(n-\sum_{i=1}^{r}|R_i|-|V^{r+1}|\right)/s$ edges in the set
	\begin{equation*}
		\bigsqcup\limits_{i=1}^r E\left(A[V_{j_1}^i\sqcup V_{j_2}^i]\right)\setminus \left(E\left(A[V_{j_1}^i]\right)\sqcup E\left(A[V_{j_2}^i]\right)\right),\quad 1\leqslant j_1 < j_2\leqslant s;
	\end{equation*}
	\item $s|V^{i+1}|$ edges between $V^{i+1}$ and $\bigsqcup\limits_{j=1}^{s} W^i_{j},$ $i\in[r].$
\end{itemize}
Now \eqref{eq:edges_init} easily follows from these computations. Lemma \ref{lemma_Asat3} is proven.
\begin{flushright}
	$\blacksquare$
\end{flushright}
\begin{lemma}\label{Hall}
	Let $C$ be a positive integer and $\mathcal{E}\subset[n]^2$ be such that, for every $i\in[n],$
	\begin{equation*}
		\left|\{j\in[n]:(i,j)\in\mathcal{E}\}\right|\leqslant C,\quad \left|\{j\in[n]:(j,i)\in\mathcal{E}\}\right|\leqslant C.
	\end{equation*}
Then the probability, that there exists a perfect matching in $G(n,n,p)$ with no edges in $\mathcal{E},$ is $1-e^{-\Omega(n)}.$
\end{lemma}
\textbf{Proof.} Define $r = n \mod 2C.$ Consider the vertex set of $G(n,n,p).$ Let $V_1,$ $V_2$ be its parts of size $n.$ Choose first $2C+r$ vertices of the first part into $U_{\floor{n/2C}}.$ By Theorem \ref{Chernoff}, whp they have at least $np^{2C+r}-\sqrt{np^{2C+r}\ln n}$ common neighbors in the second part. Expose edges adjacent to $U_{\floor{n/2C}}.$ Find their first $2C+r$ common neighbors in the second part (if exist) and assign them to a set $W_{\floor{n/2C}}.$ Observe that
\begin{equation*}
	G(n,n,p)[U_{\floor{n/2C}}\sqcup W_{\floor{n/2C}}]
\end{equation*}
is a copy of $K_{2C+r,2C+r}.$ Fix partitions $\bigsqcup_{i=1}^{\floor{n/2C}-1}U_i,$ $\bigsqcup_{i=1}^{\floor{n/2C}-1}W_i$ of $V_1\setminus U_{\floor{n/2C}}$ and $V_2\setminus W_{\floor{n/2C}}$ into sets of size $2C$ respectively.

Consider a random bipartite graph $H$ with parts $U$ and $W$ where vertices of $U$ are sets $U_i,$ $i\in[\floor{n/2C}-1],$ and vertices of $W$ are sets $W_i,$ $i\in[\floor{n/2C}-1],$ with $U_i$ adjacent to $W_j$ iff every $u\in U_i$ is adjacent to every $w\in W_j$ in $G(n,n,p).$ Then, ${\sf P}\left(U_i\sim W_j\text{ in }H\right) = p^{4C^2}.$ Denote by $B$ the event that there exists a permutation $\sigma\in S_{\floor{n/2C}-1}$ such that, for every $i\in[\floor{n/2C}-1],$ $G(n,n,p)[U_i\sqcup W_{\sigma(i)}]$ is a copy of $K_{2C,2C}.$

Let $\sigma\in S_{\floor{n/2C}-1}$ be the random permutation (say, the minimum in the lexicographical order on $S_{\floor{n/2C}-1}$) defined by $B$ if $B$ holds. Define $\sigma\in S_{\floor{n/2C}-1}$ to be identity when $B$ does not hold. Notice that $B$ holds iff $H$ has a perfect matching. By Theorem \ref{bip}, the probability that $H$ has a perfect matching is  $1-e^{-\Omega(n)}.$

Expose edges in $G(n,n,p)$ non-adjacent to $U_{\floor{n/2C}}.$ Find $\sigma\in S_{\floor{n/2C}-1}$ and the respective perfect matching in $H$ (or, in other words, the disjoint union of complete bipartite graphs $K_{2C,2C}$). 

By Hall's marriage theorem, we can find a perfect matching in each complete bipartite graph with parts $U_i,$ $W_{\sigma(i)}$ with no edges in $\mathcal{E},$ $i\in[\floor{n/2C}-1].$ Indeed, consider a bipartite graph obtained from the complete bipartite graph with parts $U_i,$ $W_{\sigma(i)}$ by removing the edges from $\mathcal{E}.$ Consider a subset $S$ of $U_i.$ If $|S|\leqslant C,$ then any vertex of $S$ has at least $C$ neighbors and, therefore, $N(U_i)\geqslant C\geqslant |U_i|.$ If $|S|>C,$ consider any vertex in $W_{\sigma(i)}.$ It has at least $C$ neighbors in $U_i.$ Hence, every vertex from $W_{\sigma(i)}$ has a neighbor in $S$ and therefore $N(S)=W_{\sigma(i)}.$

It remains to find a valid perfect matching in the bipartite graph on $U_{\floor{n/2C}}\sqcup W_{\floor{n/2C}}.$ If $r\geqslant 1,$ match the first vertex of $U_{\floor{n/2C}}$ with its arbitrary neighbor. Note that the remaning vertices still have at most $C$ non-neighbors. Repeat this $r-1$ times until $K_{2C,2C}$ remains. Proceed by Hall's marriage theorem.
\begin{flushright}
	$\blacksquare$
\end{flushright}
\begin{lemma}\label{lemma_property_m}
	Let $s$ be the minimum possible positive integer such that $(2s^2+1)(1-p)^s< 1.$ For every $p\leqslant 1-1/\sqrt[3]{7},$ whp $G(n,p)$ contains a spanning $(s,r)$-flower $A$ where 
	$$
	r=\ceil[\bigg]{\frac{5\log_{1/(1-p)}n}{8s}},\quad
	|V^2|+\ldots+|V^{r+1}|=\frac{n((1-p)^s+o(1))}{1-(1-p)^s}\quad \text{and}\quad |V^{r+1}| = o(\sqrt{n}).
	$$
\end{lemma}
\noindent\textbf{Proof.} We recursively construct sets $[n]=V^1\reflectbox{ $\subset$ }\ldots \reflectbox{ $\subset$ }V^{r}\reflectbox{ $\subset$ }V^{r+1}.$  In the same way as in the proof of Lemma \ref{lemma_property_t}, for any $s\in\mathbb{N},$ $i\in[r],$ $v_1,\ldots,v_s\in V^i$ and $\nu_k\in\{v_k,\overline{v}_k\},$ $k\in[s],$ let us denote by $N_i(\nu_1,\ldots,\nu_{s})$ the set of all common neighbors of $\{v_k,\; k\in[s]\;|\;\nu_k = v_k\}$ among common non-neighbors of $\{v_k,\; k\in[s]\;|\;\nu_k = \overline{v}_k\}$ in $V^i$ (not including $v_k \text{ such that }\nu_k = \overline{v}_k$). Let $n_i(\nu_1,\ldots,\nu_{s})$ stand for the cardinality of $N_i(\nu_1,\ldots,\nu_{s}).$

Put $\varepsilon = \frac{1-(2s^2+1)(1-p)^s}{2s}$ and let $V^1=[n].$ Assume that, for some $i\in[r],$ the sets $V^1\reflectbox{ $\subset$ }\ldots \reflectbox{ $\subset$ }V^{i}$ are already constructed and $G(n,p)[V^{i}]\stackrel{d}=G(|V^{i}|,p).$ Suppose that with probability $1-O(i/\sqrt{n})$ (the constant factor in $O(i/\sqrt{n})$ does not depend on $i$)

\begin{equation}\label{eq_small_p_1}
	|V^j|= n(1-p)^{s(j-1)}+O\left(j\sqrt{n(1-p)^{s(j-1)}\ln n}\right) \text{ for all } j\in[i].
\end{equation}
Since $i\leqslant r,$ this implies $|V^i|\gg n^{1/3}.$ All the probabilities below are conditioned on the event \eqref{eq_small_p_1}. Choose an arbitrary vertex $v_0^i$ in $V^i.$ Denote by $B_0$ the event that
$n_i(v_0^i)\geqslant |V^i|p-\sqrt{|V^i|p\ln n}.$ By Theorem \ref{Chernoff}, the probability that $B_0$ holds is $1-O(1/\sqrt{n}).$ $B_0$ (jointly with \eqref{eq_small_p_1}) implies $n_i(v_0^i)\gg n^{1/3}.$ Expose edges adjacent to $v_0^i$ in $V^i.$ Choose any $s$ neighbors of $v_0^i$ and denote them by $v_1^i,\ldots,v_s^i.$

Denote by $B_1$ the event that
$$
|V^i|p - \sqrt{|V^i|p\ln n}\leqslant n_i(v_1^i)\leqslant |V^i|p + \sqrt{|V^i|p\ln n}.
$$
By Theorem \ref{Chernoff}, the probability that $B_1$ holds is $1-O\left(1/\sqrt{n}\right).$ $B_1$ (jointly with \eqref{eq_small_p_1}) implies $n_i(v_1^i)\gg n^{1/3}$ and $n_i(\overline{v}_1^i)\gg n^{1/3}.$ Expose edges adjacent to $v_1^i$ in $V^i.$

For every $2\leqslant k\leqslant s,$ denote by $B_k$ the event that each $n_i(\nu^i_1,\ldots,\nu^i_k),$ $\nu^i_j\in\{v_j^i,\overline{v}_j^i\},$ differs from its expectation ${\sf E}n_i := {\sf E}\left(n_i(\nu_1,\ldots,\nu_k)|n_i(\nu_1,\ldots,\nu_{k-1})\right)$ by at most $\sqrt{{\sf E}n_i\ln n}.$ By Theorem \ref{Chernoff}, the probability that $B_0\cap\ldots\cap B_k$ holds is $1-O\left(1/\sqrt{n}\right).$ $B_k$ (jointly with \eqref{eq_small_p_1}) implies  $n_i(\nu^i_1,\ldots,\nu^i_k)\gg n^{1/3},$ $\nu^i_j\in\{v_j^i,\overline{v}_j^i\}.$ Expose edges adjacent to $v_k^i$ in $V^i.$

Let us construct disjoint sets $V_1^i,\ldots,V_s^i\subset V^i$ one by one in the following way. For every $j=1,\ldots,s,$ and every $T=(\nu_1^i,\nu_2^i,\ldots,v_j^i,\ldots,\nu_s^i)$ such that $\nu_k^i\in\{v_k^i,\overline{v}_k^i\},k\in[s]\setminus\{j\},$ put arbitrary $\floor{n_i(T)/\lambda(T)}$ vertices of $N_i(T)$ into $V^i_j$ (that were not exploited for the previously constructed sets $V_1^i,\ldots,V_{j-1}^i$), where $\lambda(T)$ is the number of those $\nu_j^i$ which are equal to $v_j^i$ in $T.$ Put all remaining vertices of $N_i(T)$ into $V^{i+1}.$ It follows from $B_0\cap\ldots\cap B_s$ that 
$$
|V_j^i|=|V^i|(1-(1-p)^s)/s + O\left(\sqrt{|V^i|\ln n}\right)\quad \text{for every } j\in[s],
$$
and
\begin{equation}\label{eq_ni_new}
	|V^{i+1}| = |V^i|(1-p)^s + O\left(\sqrt{|V^i|\ln n}\right).
\end{equation}

For every $j\in[s],$ in every $V_j^i$ with an odd cardinality choose any vertex, remove edges joining it to $v^i_j$ and add this vertex to $V_{i+1}.$ Now all $V_j^i$ have even cardinalities but not necessarily equal. Without loss of generality, assume that $V_1^i$ has minimum cardinality among $V_1^i,\ldots,V_s^i.$ For every $j\in\{2,\ldots,s\},$ move $|V_j^i|-|V_1^i|$ vertices from $V_j^i$ to $V^{i+1}.$ Clearly, we still have \eqref{eq_ni_new}. This condition (jointly with \eqref{eq_small_p_1}) implies $|V^{i+1}|\gg n^{1/3}.$

Remove all other edges adjacent to $v_0^i$ in $V^i$ except for $\{v_0^i,v_1^i\},\ldots,\{v_0^i,v_s^i\}.$ Remove all other edges adjacent to every $v_j^i$ except for those that connect it to $V_j^i.$

Split the sets $V^i_j,$ $j\in[s],$ into $2s$ equal parts each of size $b_i:=|V^{i+1}|$ and $2$ parts of size $y_i:=\frac{1}{2}\left(|V^i_j|-2sb_i\right).$ Denote the first $s$ parts of size $b_i$ by $U_{j,1}^i,\ldots,U_{j,s}^i.$ Put $W_j^i=U^i_{j,1}.$ Let $U_j^i$ be the union of one of the parts of size $y_i$ with $U_{j,1}^i\cup\ldots\cup U_{j,s}^i,$ and let $L_j^i=V_j^i\setminus U_j^i.$  For every $j\in[s],$ a random graph induced by $G(n,p)$ between the sets $U^i_{j}$ and $L^i_{j}$ is a bipartite random graph distributed as $G(sb_i+y_i,sb_i+y_i,p).$ By Theorem \ref{bip}, for every $j\in[s],$ the probability that it has a perfect matching is $1-e^{-\omega\left(n^{1/3}\right)}.$ For every $j\in[s],$ expose edges inside $V_j^i$ and find a perfect matching between the sets $U_j^i$ and $L_j^i$ (if exists).

For every $j_1\neq j_2\in[s],$ define 
$$
f(j_1,j_2)\in\{1,\ldots,s\},\quad f(j_1,j_2) \equiv j_2-j_1+1\mod s.
$$
The random graph induced by $G(n,p)$ between the sets $W_{j_1}^i$ and $U_{j_2,f(j_1,j_2)}^i$ is a bipartite random graph distributed as $G(b_i,b_i,p).$  Suppose that a cycle as in Figure \ref{fig:7} between $W_{j_1}^i,$ $U_{j_2,f(j_1,j_2)}^i,$ $W_{j_3}^i,$ $U_{j_4,f(j_3,j_4)}$ for distinct $j_1,j_2,j_3,j_4$ from $[s],$ is created by such bipartite graphs. It is possible only when $f(j_1,j_2) = f(j_3,j_2).$ This is equivalent to $j_1 = j_3$ and leads to a contradiction. For every $j_1\neq j_2\in[s],$ by Theorem \ref{bip}, the probability that there exists a perfect matching between the sets $W_{j_1}^i$ and $U_{j_2,f(j_1,j_2)}^i$ is $1-e^{-\omega\left(n^{1/3}\right)}.$ 

For every $j\in[s],$ a random graph induced by $G(n,p)$ between the sets $W_{j}^i$ and $V^{i+1}$ is a bipartite random graph distributed as $G(b_i,b_i,p).$ Suppose that a cycle as in Figure \ref{fig:2} between $W^i_{j_1,1},$ $U^i_{j_2,f(j_1,j_2)},$ $W^i_{j_3,1}$ for distinct $j_1,j_2,j_3$ from $[s],$ and $V^{i+1}$ is created. It is possible only when $f(j_1,j_2)=f(j_3,j_2).$ This is equivalent to $j_1=j_3$ and leads to a contradiction. For every $j\in[s],$ by Theorem \ref{bip}, the probability that there is a perfect matching between the sets $W_{j}^i$ and $V^{i+1}$ is $1-e^{-\omega\left(n^{1/3}\right)}.$

Expose edges adjacent to $W_{j}^i$ going outside $V_j^i,$ $j\in[s].$ For every $j\in[s],$ find a perfect matching between $W_{j}^i$ and $V^{i+1}$ (if exists). For every $j_1\neq j_2\in [s],$ find a perfect matching between $W_{j_1}^i$ and $U_{j_2,f(j_1,j_2)}^i$ (if exists). Remove all other edges adjacent to $W_{j}^i$ in $V^i,$ $j\in[s].$\\

Let us show by induction that, for every $j_1\neq j_2\in[s],$ the probability of the existence of a perfect matching between $L_{j_1}^i$ and $L_{j_2}^i$ so that no cycles as in Figure \ref{fig:5} are created is $1-e^{-\omega\left(n^{1/3}\right)}.$ 

Consider $L^i_1$ and $L^i_2.$ The random graph induced by $G(n,p)$ between $L^i_1$ and $L^i_2$ is a bipartite random graph distributed as $G(sb_i+y_i,sb_i+y_i,p).$ By Theorem \ref{bip}, the probability of the existence of a perfect matching in this graph is $1-e^{-\omega\left(n^{1/3}\right)}.$

Suppose that the probability of the existence of perfect matchings between all pairs of sets $L_{1}^i,\ldots,L_{\ell-1}^i,$ $3\leqslant\ell\leqslant
s,$ so that the edges of these matchings do not form any $C_4$ as in Figure \ref{fig:5} (we will call such matchings \textit{valid}), is $1-e^{-\omega\left(n^{1/3}\right)}.$ Given this event holds, let us find the probability of the existence of valid perfect matchings between the sets $L_{1}^i$ and $L_{\ell}^i,$ $L_{2}^i$ and $L_{\ell}^i,\ldots,L_{\ell-1}^i$ and $L_{\ell}^i.$ Let us do this by induction as well. Suppose that the probability of the existence of valid perfect matchings between pairs $L_{1}^i$ and $L_{\ell}^i,$ $L_{2}^i$ and $L_{\ell}^i,\ldots,L_{j-1}^i$ and $L_{\ell}^i,$ $1\leqslant j\leqslant\ell-1,$ is $1-e^{-\omega\left(n^{1/3}\right)}.$ Given this event holds, let us find the probability of the existence of a valid perfect matching between $L_{j}^i$ and $L_{\ell}^i.$

A random graph induced by $G(n,p)$ between sets $L^i_{j}$ and $L^i_{\ell}$ is distributed as a random bipartite graph $G(sb_i+y_i,sb_i+y_i,p).$ For any vertex from $L^i_{\ell},$ the cycle as in Figure \ref{fig:5} can be created iff it forms an edge with one of at most $(j-1)^3<s^3$ vertices from $L_j^i.$ By Lemma \ref{Hall}, the probability of the existence of a perfect matching in this graph without such edges is $1-e^{-\omega\left(n^{1/3}\right)}.$ This finishes the inductive argument.\\

Expose edges adjacent to $\bigsqcup\limits_{j=1}^{s}L^i_{j}.$ For all $1\leqslant j<\ell\leqslant s,$ find valid perfect matchings between $L^i_{j}$ and $L^i_{\ell}$ (if exist).\\

It remains to show by induction that, for every $j_1\neq j_2\in[s],$ the probability of the existence of a perfect matching between $U_{j_1}^i\setminus \left(U_{j_1,f(j_2,j_1)}^i\cup W_{j_1}^i\right)$ and $U_{j_2}^i\setminus \left(U_{j_2,f(j_1,j_2)}^i\cup W_{j_2}^i\right)$ so that no cycles as in Figures \ref{fig:3} and \ref{fig:7} are created, is $1-e^{-\omega\left(n^{1/3}\right)}$ (we call such matchings \textit{valid}). The proof of this fact is the same as the proof of existence of valid perfect matchings between $L$-sets. Indeed, when all matchings between pairs of sets preceeding the pair $(U_j^i,U_{\ell}^i),$ $j<\ell,$ are constructed, for any vertex from $U^i_{\ell}\setminus \left(U^i_{\ell,f(j,\ell)}\cup W_{\ell}^i\right),$ a cycle as in Figure \ref{fig:3} can be created iff it forms an edge with one particular vertex from $U^i_{j}\setminus \left(U^i_{j,f(\ell,j)}\cup W^i_{j}\right);$ a cycle as in Figure \ref{fig:7} can be created iff it forms an edge with one of at most $(j-1)^3$ vertices from $U^i_{j}\setminus \left(U^i_{j,f(\ell,j)}\cup W^i_{j}\right).$ Lemma \ref{Hall} implies the desired.\\

Expose edges between sets $U_{j_1}^i\setminus \left(U_{j_1,f(j_2,j_1)}^i\cup W_{j_1}^i\right)$ and $U_{j_2}^i\setminus \left(U_{j_2,f(j_1,j_2)}^i\cup W_{j_2}^i\right).$ For all $1\leqslant j_1<j_2\leqslant s,$ find valid perfect matchings between $U_{j_1}^i\setminus \left(U_{j_1,f(j_2,j_1)}^i\cup W_{j_1}^i\right)$ and $U_{j_2}^i\setminus \left(U_{j_2,f(j_1,j_2)}^i\cup W_{j_2}^i\right)$ (if exist). Among the exposed edges remove all other edges adjacent to $\bigsqcup\limits_{j=1}^{s}U^i_{j}$ which are not in the perfect matchings. Evidently, $G(n,p)[V^{i+1}]\stackrel{d}=G(|V^{i+1}|,p).$ Finally, it follows from \eqref{eq_small_p_1}, \eqref{eq_ni_new} that with probability $1-O\left((i+1)/\sqrt{n}\right)$
\begin{equation*}
	|V^{i+1}|= n(1-p)^{si}+O\left(i\sqrt{n(1-p)^{si}\ln n}\right)
\end{equation*}
(where the big-O is bounded by the same constant as in \eqref{eq_small_p_1}). This finishes the construction of $V^1\reflectbox{ $\subset$ }\ldots \reflectbox{ $\subset$ }V^{r}\reflectbox{ $\subset$ }V^{r+1}.$\\

Expose edges within $V^{r+1}.$ Remove some edges such that a $C_4$-saturated graph is left. It remains to estimate $|V^2|+\ldots+|V^{r+1}|$ and $|V^{r+1}|.$ It follows from \eqref{eq:edges_init} that whp
$$
|V^{r+1}| = n(1-p)^{sr}(1+o(1))=o(\sqrt{n}),
$$
$$
|V^2|+\ldots+|V^{r+1}| = \frac{n(1-p)^s}{1-(1-p)^s} + o(n).
$$
\begin{flushright}
	$\blacksquare$
\end{flushright}

Due to Lemma \ref{lemma_Asat3} and Lemma \ref{lemma_property_m}, whp $G(n,p)$ contains a subgraph $A$ with at most 
$$
n\left(\frac{s+1}{2}+\frac{s(1-p)^{s}}{1-(1-p)^s}\right)(1+o(1))
$$
edges such that $A$ is $C_4$-saturated. Hence, whp
$$
\mathrm{sat}\left(G(n,p),C_4\right)\leqslant n\left(\frac{s+1}{2}+\frac{s(1-p)^{s}}{1-(1-p)^s}\right)(1+o(1)).
$$
The inequality \eqref{cycle_4'_up} is proven.

\section{Proof of Theorem \ref{th_main_3}}\label{proof_th4}
Let us start from several auxiliary assertions.
\begin{claim}\label{lemma_c4_1}
	Let $\delta\in\left(0,p/4\right)$ be a constant, $k = k(n) > \frac{16}{p}\ln n.$ Then whp any induced subgraph of $G(n,p)$ on $k$ vertices has at least $\delta k^2$ edges.
\end{claim}
\noindent\textbf{Proof of Claim \ref{lemma_c4_1}.}
	Let $X_i$ be a random variable equal to the number of edges in the $i$-th $k$-element subset of $[n].$ It has binomial distribution $\mathrm{Bin}\left(\binom{k}{2},p\right).$ By the union bound over all possible sets of size $k$ in $[n]$ and by Theorem \ref{Chernoff}, the probability that there exists a $k$-element subset of $[n]$ that contains less than $\delta k^2$ edges is at most
	\begin{equation*}
		\sum_{i=1}^{\binom{n}{k}}{\sf P}\left(X_i<\delta k^2\right)\leqslant n^ke^{-\frac{\left(p\binom{k}{2}-\delta k^2\right)^2}{k^2p}}=\exp\left(k\ln n - k^2\frac{\left(p/2-\delta+o(1)\right)^2}{p}\right)\to 0,\quad n\to\infty.
	\end{equation*}
\begin{flushright}
	$\blacksquare$
\end{flushright}
\begin{claim}\label{claim8}
	Whp, for any set $V\subset[n]$ of size at least $\frac{3}{p}\ln n$ in $G(n,p),$ the number of vertices outside this set having at most $\frac{1}{2}\ln n$ neighbors in $V$ is less than $\ln^3n.$
\end{claim}
\noindent\textbf{Proof of Claim \ref{claim8}.} For every $v\in [n]\setminus V,$ define a random variable $X_v$ which is equal to $1$ if $v$ has at most $\frac{1}{2}\ln n$ neighbors in $V,$ and equals to $0$ otherwise. The number of neighbors of $v$ in $V$ has binomial distribution $\mathrm{Bin}\left(|V|,p\right).$ By Theorem \ref{Chernoff},
\begin{equation*}
	t:={\sf P}\left(X_v = 1\right)\leqslant e^{-\frac{\left(|V|p-\frac{1}{2}\ln n\right)^2}{2|V|p}} = o\left(\frac{1}{n}\right).
\end{equation*}
By Theorem \ref{Chernoff},
\begin{equation*}
	q:={\sf P}\left(\sum_{v\in [n]\setminus V}X_v\geqslant\ln^3 n \right) \leqslant e^{-\frac{\left(\ln^3n-(n-|V|)t\right)^2}{2(n-|V|)t+\left(\ln^3n-(n-|V|)t\right)/3}}\leqslant e^{-3\ln^3n(1+o(1))}.
\end{equation*}
By the union bound, the probability of the existence of a set $V\subset[n]$ of size at least $\frac{3}{p}\ln n$ in $G(n,p)$ and at least $\ln^3n$ vertices outside $V$ having at most $\frac{1}{2}\ln n$ neighbors in $V$ is bounded from above by
\begin{equation*}
	\binom{n}{\ceil[\big]{\frac{3}{p}\ln n}}q\leqslant e^{- 3\ln^3n(1+o(1))}\to 0,\quad n\to\infty.
\end{equation*}

\begin{flushright}
	$\blacksquare$
\end{flushright}
\begin{lemma}\label{lemma_c4_2}
	Let $c>0,$ $\varepsilon>0.$ Then whp, for any induced subgraph $H\subset G(n,p)$ on at least $\varepsilon n$ vertices and any spanning $F\subset H$ such that it is $C_4$-saturated in a spanning subgraph of $H$ obtained by removing at most $cn$ edges,
	\begin{enumerate}
		\item there are at most $\floor[\big]{\frac{3}{p}\ln n}$ isolated vertices in $F;$
		\item the set of vertices of $F$ having degree 1 has cardinality $O\left(\sqrt{|E(F)|}\ln n+\frac{n}{\ln n}\right)$.
	\end{enumerate}
\end{lemma}
\noindent\textbf{Proof of Lemma \ref{lemma_c4_2}.}
Let $F\subset H\subset G(n,p)$ be subgraphs from the condition of Lemma \ref{lemma_c4_2}. Let $W$ be the set of all isolated vertices in $F.$  Clearly, in $H,$ the number of edges having at least one vertex in $W,$ is at most $cn,$ since insertion of any of them in $F$ can not create $C_4.$ By Claim \ref{lemma_c4_1}, whp $|W|<n^{2/3}.$ Suppose that $|W|>\frac{3}{p}\ln n.$ By Claim \ref{claim8}, whp less than $\ln^3n$ vertices of $H$ have at most $\frac{1}{2}\ln n$ neighbors in $W.$ Since $|V(H)|\geqslant\varepsilon n,$ we have at least $\Omega(n\ln n)$ edges going to $W$ which leads to a contradiction. Whence the first part of Lemma \ref{lemma_c4_2} follows.\\

Let $V$ be the set of vertices of $F$ having degree 1. Let $V=V_1\sqcup\ldots\sqcup V_m\sqcup U$ be a decomposition such that $U$ is the set of end-points of all edges in $F|_V$ (i.e., vertices of the inclusion-maximum matching), $V_i$ are inclusion-maximum sets of vertices having a common neighbor outside $V$. Let $\tilde H\subset H$ be obtained from $H$ by removing at most $cn$ edges in a way such that $F$ is $C_4$-saturated in $\tilde H$. Clearly, every set $V_i$ induces an empty graph in $\tilde H$. Also, there are no edges with end-points in $\tilde H|_U$ other than those that are from the matching induced by $U$ in $F$. Therefore, the number of edges in $H|_U$ is at most $cn+|U|/2$. By Claim \ref{lemma_c4_1}, whp either $|U|\leq\frac{16}{p}\ln n$ or $\frac{p}{8}|U|^2<cn+|U|/2$. Therefore, whp $|U|=O(\sqrt{n})$. \\

Let $v=|V_1|+\ldots+|V_m|=|V|-|U|$. 

For every $i\in[m]$, choose an arbitrary $v_i\in V_i$ and consider $\tilde V=\{v_1,\ldots,v_m\}$. For any edge $\{v_i,v_j\}$ from $\tilde H|_{\tilde V}$ there is a unique edge in $F$ that recovers $\{v_i,v_j\}$. Moreover, there is one-to-one correspondence between edges from $\tilde H|_{\tilde V}$ and the edges from $F$ that recover them. Therefore, $|E(F)|\geq|E(\tilde H|_{\tilde V})|$. By Claim \ref{lemma_c4_1}, whp 
\begin{equation}
	\text{either}\quad m\leq\frac{16}{p}\ln n,\quad\quad\text{or}\quad\frac{p}{8}m^2-cn\leq|E(F)|.
	\label{m_bound}
\end{equation}
Notice that, if $v\leq \frac{16m}{p}\ln n$, then 
$$
\text{ either}\quad v\leq \frac{256}{p^2}\ln^2n,\quad\quad\text{or}\quad v\leq \frac{16}{p}\sqrt{\frac{8}{p}\left(|E(F)|+cn\right)}\ln n,
$$ 
and we are done.\\

Assume that $v>\frac{16m}{p}\ln n$. Let $\ell$ be the number of sets from $V_1,\ldots,V_m$ that have cardinality at most $\frac{16}{p}\ln n$. Without loss of generality, assume that these sets are $V_{m-\ell+1},\ldots,V_m$. Due to the restriction on $v,$ we have $\ell<m$. By Claim \ref{lemma_c4_1}, whp every $H|_{V_i}$, $i\in[m-\ell]$, has at least $\frac{p}{8}|V_i|^2$ edges. Then $\frac{p}{8}\sum_{i=1}^{m-\ell}|V_i|^2\leq cn$. Notice that $v>\frac{16\ell}{p}\ln n$. Since 
$$
\sum_{i=1}^{m-\ell}|V_i|^2\geq(m-\ell)\left(\frac{\sum_{i=1}^{m-\ell}|V_i|}{m-\ell}\right)^2\geq\frac{\left(v-\ell\frac{16}{p}\ln n\right)^2}{m-\ell},
$$
we get
\begin{equation}
	\frac{\left(v-\ell\frac{16}{p}\ln n\right)^2}{m-\ell}\leq \frac{8}{p}cn.
	\label{g_ell}
\end{equation}
Letting $g(x)=\frac{\left(v-x\frac{16}{p}\ln n\right)^2}{m-x}$, we get $\ln g=2\ln\left(v-x\frac{16}{p}\ln n\right)-\ln(m-x)$ and 
$$
(\ln g)^{\prime}=-\frac{32\ln n}{p\left(v-x\frac{16}{p}\ln n\right)}+\frac{1}{m-x}=\frac{pv+16x\ln n-32m\ln n}{p\left(v-x\frac{16}{p}\ln n\right)(m-x)}.
$$
Therefore, $g$ achieves its minimum at $x=2m-\frac{pv}{16\ln n}$ (notice that both $v-x\frac{16}{p}\ln n$ and $m-x$ are positive). Thus, from \eqref{g_ell},
$$
\frac{4\left(v-m\frac{16}{p}\ln n\right)^2}{\frac{pv}{16\ln n}-m}\leq g(\ell)\leq \frac{8}{p}cn.
$$
Finally, $v\leq\frac{16}{p}m\ln n+\frac{c}{8}\frac{n}{\ln n}$. Together with \eqref{m_bound}, it gives the desired bound on $v$.
\begin{flushright}
	$\blacksquare$
\end{flushright}

\begin{lemma}\label{lemma_c4_3}
	Let $\varepsilon>0.$ Then whp, for any induced subgraph $H\subset G(n,p)$ on at least $\varepsilon n$ vertices and any spanning $F\subset H$ such that it is $C_4$-saturated in a spanning subgraph of $H$ obtained by removing at most $n$ edges,
	\begin{enumerate}
		\item $F$ does not have induced $P_s,$ $s>\frac{3}{p}\ln n + 6;$
		\item the number of induced inclusion-maximal $P_s,$ $5\leqslant s \leqslant\frac{3}{p}\ln n + 6,$ is at most $\sqrt{5n/p}.$
	\end{enumerate}
\end{lemma}
\noindent\textbf{Proof of Lemma \ref{lemma_c4_3}.}
	Assume that there exists $P_s,$ $s=\ceil[\big]{\frac{3}{p}\ln n} + 6.$ Let $U$ be the set of its $s-6\geqslant\frac{3}{p}\ln n$ central vertices (all but those that are at distance at most $2$ from the ends of the path). By Claim \ref{claim8}, whp at most $\ln^3n$ vertices outside $U$ have at most $\frac{1}{2}\ln n$ neighbors in $U.$ Therefore, whp, the number of edges between $V(H)\setminus U$ and $U$ in $H$ is $\Omega\left(n\ln n\right)$ which is bigger than $n.$ This leads to a contradiction since none of these edges can be recovered from $F.$
	
\begin{figure}
	\centering
	\includegraphics[width=0.6\linewidth]{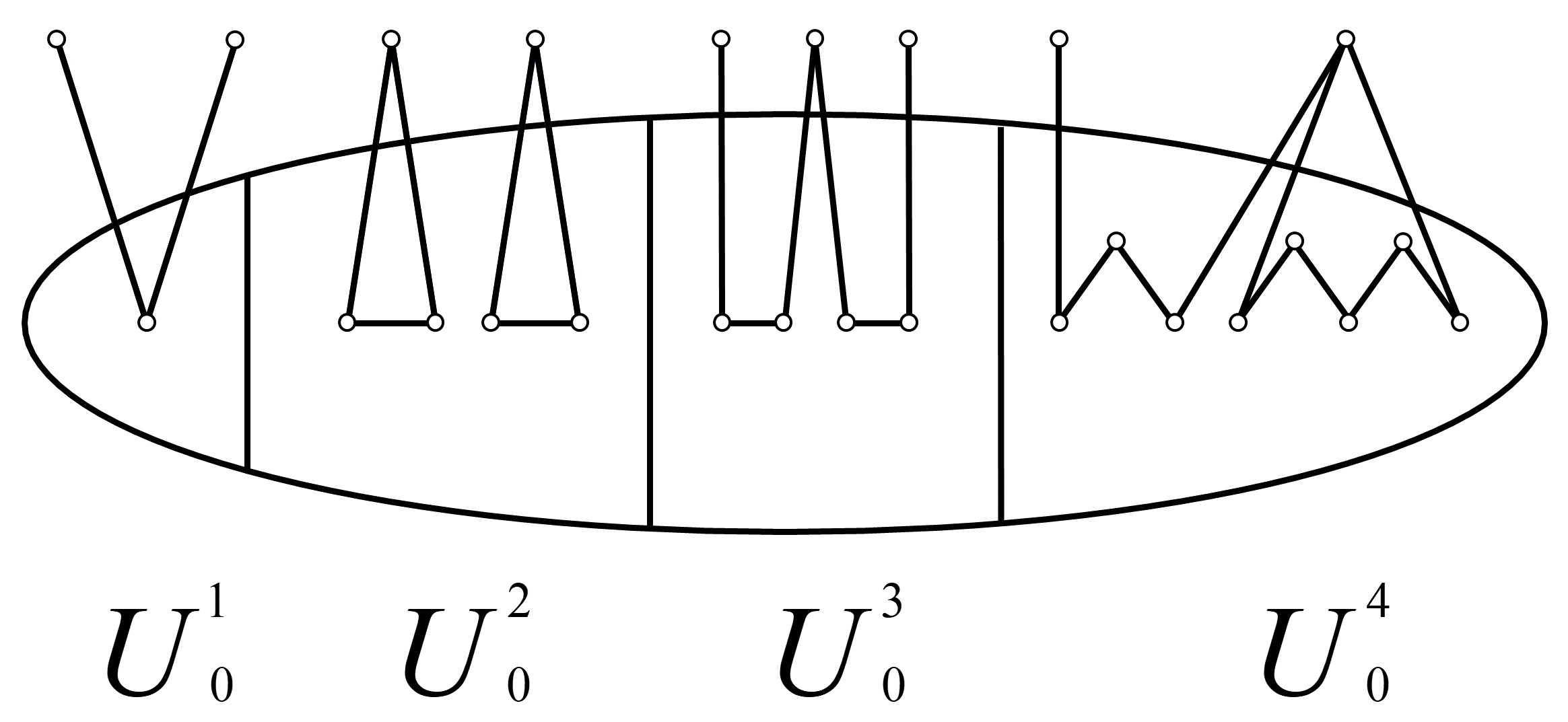}
	\caption{}
	\label{fig:U0}
\end{figure}

	Now consider the set of all inclusion-maximal $P_s,$ $s\geqslant 5.$ Select a central vertex from every such path into a set $W.$ Any edge from $H\vert_W$ can not be recovered from $F.$ If $|W|>\sqrt{5n/p},$ then, by Claim \ref{lemma_c4_1}, the number of edges in $H\vert_W$ is bigger than $n.$ We come into a contradiction.
\begin{flushright}
	$\blacksquare$
\end{flushright}
Now, let $H$ be $C_4$-saturated in $G(n,p).$ Assume that $|E(H)|\leqslant 3n/2,$ $\varepsilon>0$ is small enough. Let $U_0$ be the set of vertices having degree $2$ in $H.$ Let $U_0=U_0^1\sqcup U_0^2\sqcup U_0^3\sqcup U_0^4$ be a partition (see Figure \ref{fig:U0}), where
\begin{itemize}
	\item $U_0^1$ is the set of isolated vertices in $H\vert_{U_0},$
	\item $U_0^2$ contains all pairs of adjacent vertices that have a common neighbor outside $U_0,$
	\item $U_0^3$ contains all pairs of adjacent vertices that have different neighbors outside $U_0,$
	\item $U_0^4$ contains all the other vertices of $U_0.$
\end{itemize}
Assume that $|U_0^2\sqcup U_0^3|>\varepsilon n.$ Remove all vertices of $U_0^2$ and edges adjacent to them from $H$ and $G(n,p).$ Let us denote the obtained graphs by $\hat{H}$ and $\hat{G}$ respectively. Clearly, $\hat{H}$ is $C_4$-saturated in $\hat{G}.$ Remove all vertices of $U_0^3$ and edges adjacent to them from $\hat{H}$ and $\hat{G}.$ Let us denote the obtained graphs by $\tilde{H}$ and $\tilde{G}$ respectively. Observe that at most $\frac{1}{2}|U_0^3|$ edges should be removed from $\tilde{G}$ to make $\tilde{H}$ a $C_4$-saturated graph in $\tilde{G}.$ Set $H_1 = \tilde{H},$ $G_1 = \tilde{G}.$ Define the sets $U_0^2, U_0^3$ for these graphs. If $|U_0^2\sqcup U_0^3|>\varepsilon n,$ then do the same deletions as above and construct $H_2, G_2.$ After $m \leqslant \floor{1/\varepsilon}+1$ such steps we obtain $H_m$ and $G_m$ with $|U_0^2\sqcup U_0^3|\leqslant\varepsilon n.$ At most $\frac{1}{2}n$ edges should be removed from $G_m$ to make $H_m$ a $C_4$-saturated graph in $G_m.$ Notice that the number of deleted edges is $1.5$ times greater than the number of deleted vertices. Therefore, we may assume that $|V(H_m)|\geqslant \sqrt{\varepsilon}n.$ Otherwise, $|E(H)|>\left(\frac{3}{2} - \frac{3}{2}\sqrt{\varepsilon}\right)n,$ and there is nothing to prove.

Consider the new set $U_0$ (defined for $H_m$) and the new (defined as above) partition $U_0=U_0^1\sqcup U_0^2\sqcup U_0^3\sqcup U_0^4.$ Assume that $|U_0^1|\geqslant \varepsilon n.$ Otherwise, by Lemma \ref{lemma_c4_2} and Lemma \ref{lemma_c4_3}, we have whp
\begin{multline*}
	|E(H_m)|>\frac{3}{2}\left(|V(H_m)|-|U_0|-o(n)\right)+|U_0|+o(n) = \\
	 \frac{3}{2}\;|V(H_m)|-\varepsilon n - \frac{1}{2}\;|U_0^4|+o(n) = \frac{3}{2}\;|V(H_m)|-\varepsilon n + o(n)
\end{multline*}
implying that $|E(H)|>\frac{3}{2}\;n - \varepsilon n + o(n).$

Let $C = \floor{1/2\varepsilon}.$ Let $W:=N(U_0^1)$ be the set of all neighbors of vertices from $U_0^1$ in $H_m.$ For $v\in W,$ denote by $\mathrm{deg}_{*} v$ the number of neighbors of $v$ in $U_0^1$ and denote by $\mathrm{deg}^{*} v$ the number of neighbors of $v$ not in $U_0^1.$ Let $W = W_{*}\sqcup W^{*},$ where $W_{*}$ is the set of all $v$ with $\mathrm{deg}_{*} v\leqslant C.$ Let $U_{*}\subset U_0^1$ be the set of all vertices with both neighbors in $W_{*}.$ 
\begin{figure}
	\centering
	\includegraphics[width=0.6\linewidth]{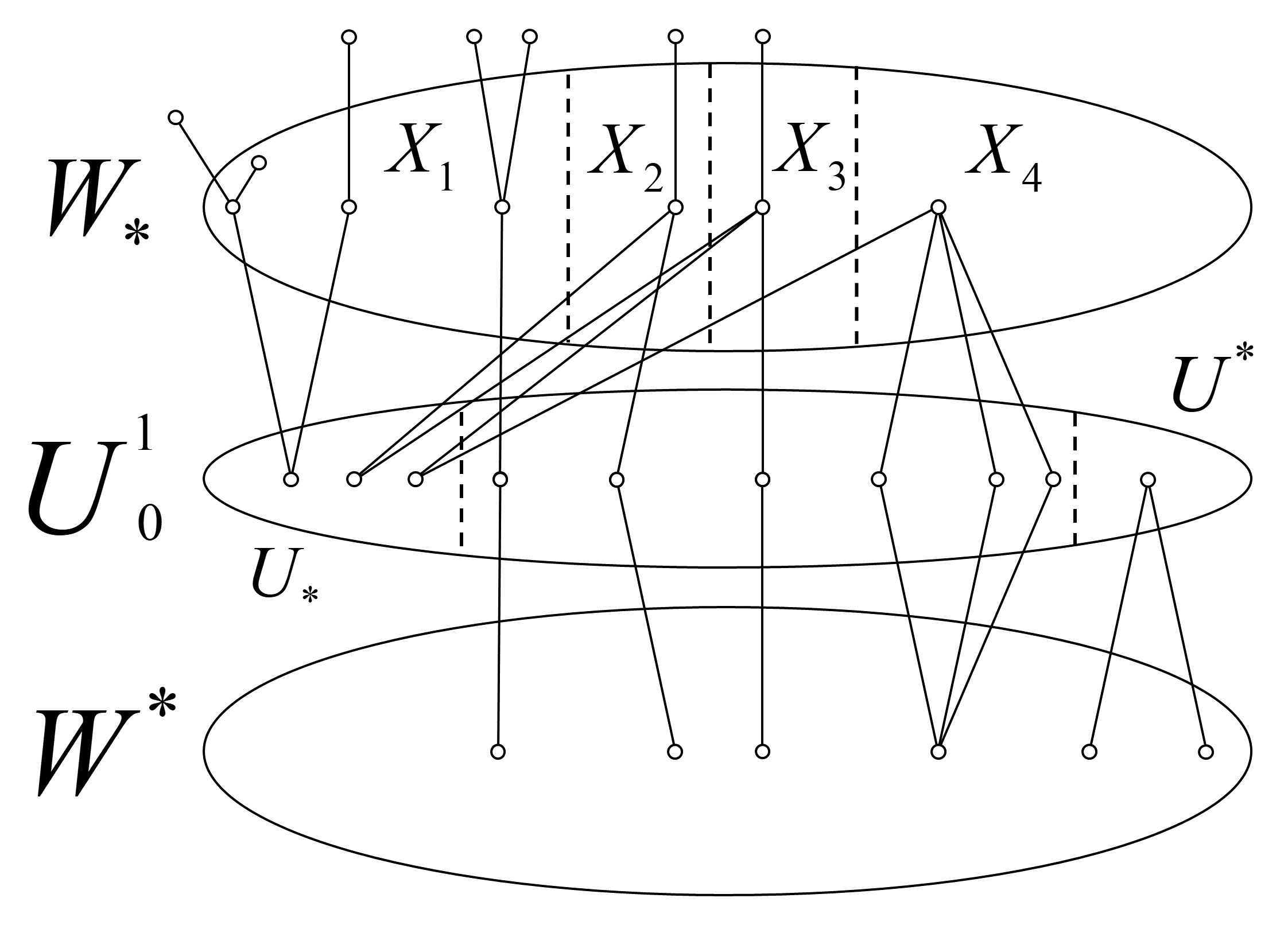}
	\caption{}
	\label{fig:W}
\end{figure}
Any edge $\{v_1,v_2\},$ $v_1,v_2\in W_{*},$ may recover at most $C^2$ edges from $G_m\vert_{U_{*}}.$ Denote $|U_{*}|$ by $k.$ Let us prove that whp $k\leqslant 4C\sqrt{n/p}.$ Assume that $k>\frac{16}{p}\ln n.$ By Claim \ref{lemma_c4_1}, there are at least $pk^2/8$ edges in $G\vert_{U_{*}}.$ Therefore, the number of edges in $
H_m\vert_{W_{*}}$ is at least $(pk^2/8-n)/C^2.$ Then $(pk^2/8 - n)/C^2\leqslant 3n/2,$ and so $k\leqslant 4C\sqrt{n/p}$ as desired. It is also obvious that $C|W^{*}|\leqslant 2|U_0^1|.$ Then, $|W^{*}|\leqslant 2n/C.$ 

If $u\in U_0^1\setminus U_{*},$ then either it has both neighbors in $W^{*},$ or a unique neighbor $v\in W_{*}.$ In the latter case, $\mathrm{deg}\, v\geqslant 3$ (since $v\notin U_0$). Therefore, either $\mathrm{deg}_{*}v=1$ and $\mathrm{deg}^{*}v\geqslant 2,$ or $\mathrm{deg}_{*}v=2$ and $\mathrm{deg}^{*}v\geqslant 1,$ or $\mathrm{deg}_{*}v\geqslant 3.$ For $i\in\{1,\ldots,C\},$ let $X_i$ be the set of vertices $v$ in $W_{*}$ with $\mathrm{deg}_{*}v = i$ (see Figure \ref{fig:W}). Let $U^{*}$ be the set of vertices from $U_0^1$ having both neighbors in $W^{*},$ $z=|U^{*}|.$ Let $I$ be the set of vertices with degree at most one in $H_m.$ Let $r$ be the number of vertices of $H_m$ outside $U_0\sqcup W\sqcup I.$ Let us bound from below $|E(H_m)|.$ Since
\begin{itemize}
	\item every vertex from $X_1\sqcup X_2$ has degree at least $3,$
	\item every vertex from $X_i,$ $i\geqslant 3,$ has degree at least $i,$
	\item vertices from $U_0^1$ receive $\sum_{i=1}^{C}ix_i$ edges from $X_1\sqcup\ldots\sqcup X_C,$
	\item vertices from $U_0^1\setminus\left(U_{*}\sqcup U^{*}\right)$ send the same number of edges to $W^{*}$ as they receive from $W_{*}$ (this number equals $\sum_{i=1}^{C}ix_i -2k$),
	\item vertices from $U^{*}$ send $2z$ edges to the vertices from $W^{*},$
	\item vertices outside $U_0\sqcup W\sqcup I$ have degrees at least $3,$
\end{itemize}
we have
\begin{equation*}
	|E(H_m)|\geqslant\frac{1}{2}\left(3|X_1|+3|X_2|+\sum_{i=3}^{C}i|X_i|+\sum_{i=1}^Ci|X_i|+3r\right)+\sum_{i=1}^Ci|X_i|-2k+2z.
\end{equation*}
Moreover,
\begin{equation*}
	|V(H_m)|\leqslant \sum_{i=1}^C|X_i|+\sum_{i=1}^Ci|X_i|+z+|W^{*}|+|U_0^2|+|U_0^3|+|U_0^4|+|I|+r.
\end{equation*}

Then, by Lemma \ref{lemma_c4_2} and Lemma \ref{lemma_c4_3},
\begin{multline*}
	\frac{|E(H_m)|}{|V(H_m)|}\geqslant \frac{3|X_1|+9|X_2|/2+2\sum_{i=3}^Ci|X_i|+3r/2+2z-o(n)}{2|X_1|+3|X_2|+\sum_{i=3}^C(i+1)|X_i|+r+z+n\left(2/C+\varepsilon+o(1)\right)} \geqslant\\ \frac{3}{2} - \frac{n\left(2/C+\varepsilon+o(1)\right)}{|V(H_m)|}\geqslant \frac{3}{2}-\frac{2/C+\varepsilon}{\sqrt{\varepsilon}}-o(1)\geqslant\frac{3}{2}-6\sqrt{\varepsilon}-o(1)
\end{multline*}
implying that $|E(H)|\geqslant \frac{3}{2}n-(6\sqrt{\varepsilon}-o(1))|V(H_m)|\geqslant \frac{3}{2}n-6\sqrt{\varepsilon}n+o(n).$
\section{Acknowledgements}
Yury Demidovich is supported by the Ministry of Science and Higher Education of the Russian Federation in the framework of MegaGrant no 075-15-2019-1926 and by the Simons Foundation. Maksim Zhukovskii is supported by RFBR and INSF, grant number 20-51-56017.

\end{document}